\documentclass[12pt,reqno]{amsart}

\usepackage{listings} 
\usepackage{xcolor} 

\usepackage{hyperref}
\usepackage{amsthm}
\usepackage[letterpaper]{geometry}
\usepackage{geometry}

\newtheorem{theorem}{Theorem}[section]
\newtheorem{proposition}[theorem]{Proposition}

\theoremstyle{remark}

\theoremstyle{definition}
\newtheorem{definition}[theorem]{Definition}

\newcommand{\image}{\mathrm{image}}
\newcommand{\R}{\mathbb{R}} 
\newcommand{\Z}{\mathbb{Z}}

\newcommand{\N}{\mathbb{N}}
\newcommand{\p}{\partial}
\newcommand{\Con}{\text{Con}}
\newcommand{\Vol}{\text{Vol}}
\allowdisplaybreaks[4]

\title{Harmonic forms on ALE Ricci-flat 4-manifolds}

\author{Gao Chen} \thanks{The first named author is supported by the Project of Stable Support for Youth Team in Basic Research Field, Chinese Academy of Sciences, YSBR-001.}
\address{Institute of Geometry and Physics, University of Science and Technology of China, 96 Jinzhai Road, Baohe District, Hefei, China}
\email{chengao1@ustc.edu.cn}
\author{Hao Yan}
\address{Department of Mathematics, University of California, Irvine, Irvine, CA 92697, USA}
\email{hyan27@uci.edu}

\begin{document}
	
	\begin{abstract}
		In this paper, we compute the expansion of some harmonic functions and 1-forms on ALE Ricci-flat 4-manifolds.
	\end{abstract}
	
	\maketitle
	
	\section{Introduction}
	
	This paper studies Asymptotic Locally Euclidean (ALE) Ricci-flat 4-manifolds (see Def \ref{defALE}). In \cite{BaKaNa}, Bando, Kasue, and Nakajima studied the expansion of the metric near infinity and conjectured that all simply-connected ALE Ricci-flat 4-manifolds $(X,g_X)$ must be hyper-K\"ahler, which have been classified by Kronheimer (see \cite{K1} and \cite{K2}). Progress on this conjecture has been limited, with notable contributions in \cite{LV} and \cite{M.Li}. More recently, in \cite{BH}, Biquard and Hein improved the expansion of the metric near infinity and proved that the renormalized volume $\mathcal{V}=\lim_{R\to\infty} \Vol(B_R, g_X)-\Vol(B_R, g_{\mathbb{R}^4/\Gamma})$ is non-positive. They utilized a function with a constant Laplacian to establish their results.\\
	
	Our first theorem generalized their results to more functions.
	
	\begin{theorem}
		\label{theoremfunc}
		We fix $\epsilon\in(0,1)$. For any $a_{ij}\in \mathbb{R}$, $i,j=1,\cdots,4$, such that $a_{ij}=a_{ji}$ and $\sum_{i,j=1}^{4} a_{ij}x^i x^j$ is invariant under $\Gamma$, there exists a unique smooth function $u_a$ on $X$ such that
		\begin{align}
			\Pi^*\Phi^* u_a-\sum_{i,j=1}^{4} a_{ij}x^i x^j\in W^{k,2}_{-2+\epsilon}(\R^4\backslash B_{2R}(0))
		\end{align}
		for all $k\ge 0$ (see Def \ref{defnorm}), and $\Delta_X u_a =(d\delta_X+\delta_X d)u_a= -2\sum_{i,j=1}^{4} a_{ij}\delta_{ij}$ on $X$. Moreover, the expansion of $u_a$ is given by the following:
		\begin{align}
			\Pi^*\Phi^*u_a-\sum_{i,j=1}^{4} a_{ij}\big(x^i x^j-\tilde{\eta}_{ij}-\frac{|\Gamma|\mathcal{V}\delta_{ij}}{2\pi^2r^2}\big)\in W^{k,2}_{-3+\epsilon}(\R^4\backslash B_{2R}(0)),
		\end{align}
		where $\tilde{\eta}_{ij}$ is defined in \eqref{defetaii} and \eqref{defetaij}.
	\end{theorem}	
	
	In \cite{BH}, Biquard and Hein also sketched the computation of expansions of all harmonic 1-forms asymptotic to $x^idx^j-x^jdx^i$, and used this to study Killing fields. Another main goal of this paper is to generalize their result to harmonic 1-forms asymptotic to $x^idx^j$.
	
	\begin{theorem}
		\label{theorem1form}
		We fix $\epsilon\in(0,1)$. For any $a_{ij}\in \mathbb{R}$, $i,j=1,\cdots,4$, such that $\sum_{i,j=1}^{4} a_{ij}x^i dx^j$ is invariant under $\Gamma$, there exists a unique smooth 1-form $\omega_a$ on $X$ such that
		\begin{align}
			\Pi^*\Phi^* \omega_a-\sum_{i,j=1}^{4} a_{ij}x^i dx^j\in W^{k,2}_{-3+\epsilon}(\R^4\backslash B_{2R}(0))
		\end{align}
		for all $k\ge 0$, and $\Delta_X \omega_a =(d\delta_X+\delta_X d)\omega_a=0$ on $X$. Moreover, the expansion of $\omega_a$ is given by the following:
		\begin{align}
			\Pi^*\Phi^*\omega_a-\sum_{i,j=1}^{4} a_{ij}\big(x^i dx^j-\tilde{\mu}_{ij}-\sum_{k,l=1}^4\frac{\Con[i,j,k,l]x^k}{r^4}dx^l\big)\in W^{k,2}_{-4+\epsilon}(\R^4\backslash B_{2R}(0)),
		\end{align}
		where $\tilde{\mu}_{ij}$ is defined in \eqref{defmuij}, and $\Con[i,j,k,l]$ satisfies the equations in Section 5.
	\end{theorem}	
	
	\subsection{Acknowledgements} The authors would like to thank Xiuxiong Chen, Hans-Joachim Hein, Mingyang Li and Bing Wang for numerous helpful discussions and Peirong Luo for his assistance in Wolfram Mathematica programming.
	
	\section{The settings of ALE-manifolds}
	
	In this section, we prepare the basic elements for our calculation.
	\begin{definition}
		\label{defALE}
		Let $(X,g)$ be a Ricci-flat 4-manifold. We call it ALE, if there exists a finite subgroup $\Gamma$ of $SO(4)$ acting freely on $\mathbb{S}^3$, together with a quotient map $\Pi:\R^4\backslash B_R(0)\to(\mathbb{R}^4\backslash B_R(0))/\Gamma$, and a diffeomorphism $\Phi:(\mathbb{R}^4\backslash B_R(0))/\Gamma\to X\backslash U$ for some bounded open subset $U\subset X$ such that for all $k\in \mathbb{N}_0$,
		\begin{align}
			|\nabla^k_{g_0}(\Pi^*\Phi^*g-g_0)|_{g_0}=O(r^{-4-k}) \text{ as } r\to\infty,    
		\end{align}
		
		where $x^i\ (i=1,\cdots,4)$ denotes the coordinate functions on $\R^4$, $r=\sqrt{\sum_{i=1}^4(x^i)^2}$, and $g_0$ denotes the Euclidean metric on $\mathbb{R}^4$.
	\end{definition}
	By Theorem B in \cite{BH}, for any $k_0\in \N$, we can modify $\Phi$ such that there exists a decomposition
	
	\begin{align}
		\label{expansionmetric}
		\Pi^*\Phi^*g-g_0=h+h',
	\end{align}
	where the leading term
	\begin{align}
		h:\R^4\backslash B_R(0)\to \text{Sym}^2\R^4
	\end{align}
	is a $\Gamma$-equivariant harmonic tensor that decays at infinity, and
	\begin{align}
		h=h^+ + h^-;\ \ \sum_{k=0}^{k_0}r^k|\nabla^k_{g_0}h'|_{g_0}\leq C(k_0)r^{-5}.
	\end{align}
	In the formula above, $h^+$ is a symmetric 2-tensor on $\R^4\backslash B_R(0)$ of the form
	\begin{align}
		-\frac{3}{2}r^6 h^+ &=\zeta_{11}(2\alpha^2_1-\alpha^2_2-\alpha^2_3)+\zeta_{22}(2\alpha^2_2-\alpha^2_3-\alpha^2_1)+\zeta_{33}(2\alpha^2_3-
		\alpha^2_1-\alpha^2_2)\\
		&+\zeta_{12}(\alpha_1\odot\alpha_2)+\zeta_{13}(\alpha_1\odot\alpha_3)+\zeta_{23}(\alpha_2\odot\alpha_3),
	\end{align}
	where $f\odot g=f\otimes g+g\otimes f$ denotes the symmetric product, $\alpha_j=I_j^*(rdr)$, and $(I_1,I_2,I_3)$ is the standard triple of complex structures on $\R^4$ given by the following:
	\begin{align}
		I_1^*(dx^1,dx^2,dx^3,dx^4)&=(-dx^2,dx^1,-dx^4,dx^3),\\
		I_2^*(dx^1,dx^2,dx^3,dx^4)&=(-dx^3,dx^4,dx^1,-dx^2),\\
		I_3^*&=I_2^*I_1^*=(I_1I_2)^*.
	\end{align}
	Here $(\zeta_{ij})$ is a symmetric $3\times 3$ matrix, and $h^-=\mathcal{R}^*h^+$ for some $h^+$ of the form above and some $\mathcal{R}\in\text{O(4)$\backslash$SO(4)}$ acting on $T_p\R^4=\R^4$.\\
	Let us choose $\mathcal{R}=\text{diag}(1,-1,-1,-1)$. We will use $\xi$ to denote the symmetric $3\times 3$ matrix of coefficients appearing in $h^-=\mathcal{R}^*h^+$. Thus we have the following explicit expression of $h^+$ and $h^-$ as a matrix with respect to the basis $\{dx^i\otimes dx^j\}$ in the region $\{r\geq 2R\}$:
	\begin{align}
		(-\frac{3}{2}r^6h^+)_{11}&=(x^2)^2 (2 \zeta_{11}-\zeta_{22}-\zeta_{33})+6 x^2 (\zeta_{12} x^3-\zeta_{13} x^4)\\
		&-(x^3)^2 (\zeta_{11}-2 \zeta_{22}+\zeta_{33})-6 \zeta_{23} x^3 x^4-(x^4)^2 (\zeta_{11}+\zeta_{22}-2 \zeta_{33}),\\
		(-\frac{3}{2}r^6h^+)_{12}&=x^1 (x^2 (-2 \zeta_{11}+\zeta_{22}+\zeta_{33})-3 \zeta_{12} x^3+3 \zeta_{13} x^4)\\
		&-3 \left(\zeta_{13} x^2 x^3+\zeta_{12} x^2 x^4+\zeta_{23} (x^3)^2+\zeta_{22} x^3 x^4-\zeta_{33} x^3 x^4-\zeta_{23} (x^4)^2\right),\\
		(-\frac{3}{2}r^6h^+)_{13}&=x^1 (-3 \zeta_{12} x^2+x^3 (\zeta_{11}-2 \zeta_{22}+\zeta_{33})+3 \zeta_{23} x^4)\\
		&+3 \left(\zeta_{13} (x^2)^2+\zeta_{23} x^2 x^3+x^2 x^4 (\zeta_{11}-\zeta_{33})+x^4 (\zeta_{12} x^3-\zeta_{13} x^4)\right),\\
		(-\frac{3}{2}r^6h^+)_{14}&=3 x^2 (\zeta_{13} x^1+x^3 (\zeta_{22}-\zeta_{11})-\zeta_{23} x^4)+3 \zeta_{23} x^1 x^3+x^1 x^4 (\zeta_{11}+\zeta_{22}-2 \zeta_{33})\\
		&+3 \zeta_{12} (x^2)^2+3 x^3 (\zeta_{13} x^4-\zeta_{12} x^3),\\
		(-\frac{3}{2}r^6h^+)_{22}&=(x^1)^2 (2 \zeta_{11}-\zeta_{22}-\zeta_{33})+6 x^1 (\zeta_{13} x^3+\zeta_{12} x^4)-(x^3)^2 (\zeta_{11}+\zeta_{22}-2 \zeta_{33})\\
		&+6 \zeta_{23} x^3 x^4-(x^4)^2 (\zeta_{11}-2 \zeta_{22}+\zeta_{33}),\\
		(-\frac{3}{2}r^6h^+)_{23}&=3 \zeta_{12} (x^1)^2-3 x^2 (\zeta_{13} x^1+\zeta_{23} x^4)+3 x^1 (\zeta_{23} x^3+x^4 (\zeta_{22}-\zeta_{11}))\\
		&+x^2 x^3 (\zeta_{11}+\zeta_{22}-2 \zeta_{33})-3 x^4 (\zeta_{13} x^3+\zeta_{12} x^4),\\
		(-\frac{3}{2}r^6h^+)_{24}&=-3 \zeta_{13} (x^1)^2-3 x^1 (\zeta_{12} x^2+x^3 (\zeta_{33}-\zeta_{11})+\zeta_{23} x^4)-3 \zeta_{23} x^2 x^3\\
		&+x^2 x^4 (\zeta_{11}-2 \zeta_{22}+\zeta_{33})+3 x^3 (\zeta_{13} x^3+\zeta_{12} x^4),\\
		(-\frac{3}{2}r^6h^+)_{33}&=-\left((x^1)^2 (\zeta_{11}-2 \zeta_{22}+\zeta_{33})\right)-6 x^1 (\zeta_{23} x^2+\zeta_{12} x^4)\\
		&-(x^2)^2 (\zeta_{11}+\zeta_{22}-2 \zeta_{33})+6 \zeta_{13} x^2 x^4+(x^4)^2 (2 \zeta_{11}-\zeta_{22}-\zeta_{33}),\\
		(-\frac{3}{2}r^6h^+)_{34}&=3 \left(x^1 x^2 (\zeta_{33}-\zeta_{22})+x^1 (\zeta_{13} x^4-\zeta_{23} x^1)+\zeta_{23} (x^2)^2+\zeta_{12} x^2 x^4\right)\\
		&+x^3 (3 \zeta_{12} x^1-3 \zeta_{13} x^2+x^4 (-2 \zeta_{11}+\zeta_{22}+\zeta_{33})),\\
		(-\frac{3}{2}r^6h^+)_{44}&=-\left((x^1)^2 (\zeta_{11}+\zeta_{22}-2 \zeta_{33})\right)+x^1 (6 \zeta_{23} x^2-6 \zeta_{13} x^3)\\
		&-(x^2)^2 (\zeta_{11}-2 \zeta_{22}+\zeta_{33})-6 \zeta_{12} x^2 x^3+(x^3)^2 (2 \zeta_{11}-\zeta_{22}-\zeta_{33}),\\
		(-\frac{3}{2}r^6h^-)_{11}&=(x^2)^2 (2 \xi_{11}-\xi_{22}-\xi_{33})+6 x^2 (\xi_{12} x^3-\xi_{13} x^4)-(x^3)^2 (\xi_{11}-2 \xi_{22}+\xi_{33})\\
		&-6 \xi_{23} x^3 x^4-(x^4)^2 (\xi_{11}+\xi_{22}-2 \xi_{33}),\\
		(-\frac{3}{2}r^6h^-)_{12}&=x^1 (x^2 (-2 \xi_{11}+\xi_{22}+\xi_{33})-3 \xi_{12} x^3+3 \xi_{13} x^4)\\
		&+3 \left(\xi_{13} x^2 x^3+\xi_{12} x^2 x^4+\xi_{23} (x^3)^2+\xi_{22} x^3 x^4-\xi_{33} x^3 x^4-\xi_{23} (x^4)^2\right),\\
		(-\frac{3}{2}r^6h^-)_{13}&=x^1 (-3 \xi_{12} x^2+x^3 (\xi_{11}-2 \xi_{22}+\xi_{33})+3 \xi_{23} x^4)\\
		&-3 \left(\xi_{13} (x^2)^2+\xi_{23} x^2 x^3+x^2 x^4 (\xi_{11}-\xi_{33})+x^4 (\xi_{12} x^3-\xi_{13} x^4)\right),\\
		(-\frac{3}{2}r^6h^-)_{14}&=3 x^2 (\xi_{13} x^1+x^3 (\xi_{11}-\xi_{22})+\xi_{23} x^4)+3 \xi_{23} x^1 x^3+x^1 x^4 (\xi_{11}+\xi_{22}-2 \xi_{33})\\
		&-3 \xi_{12} (x^2)^2+3 x^3 (\xi_{12} x^3-\xi_{13} x^4),\\
		(-\frac{3}{2}r^6h^-)_{22}&=(x^1)^2 (2 \xi_{11}-\xi_{22}-\xi_{33})-6 x^1 (\xi_{13} x^3+\xi_{12} x^4)-(x^3)^2 (\xi_{11}+\xi_{22}-2 \xi_{33})\\
		&+6 \xi_{23} x^3 x^4-(x^4)^2 (\xi_{11}-2 \xi_{22}+\xi_{33}),\\
		(-\frac{3}{2}r^6h^-)_{23}&=3 \xi_{12} (x^1)^2+3 \xi_{13} x^1 x^2-3 \xi_{23} x^1 x^3+3 x^1 x^4 (\xi_{11}-\xi_{22})\\
		&+x^2 x^3 (\xi_{11}+\xi_{22}-2 \xi_{33})-3 \xi_{23} x^2 x^4-3 x^4 (\xi_{13} x^3+\xi_{12} x^4),\\
		(-\frac{3}{2}r^6h^-)_{24}&=-3 \xi_{13} (x^1)^2+3 x^1 (\xi_{12} x^2+x^3 (\xi_{33}-\xi_{11})+\xi_{23} x^4)-3 \xi_{23} x^2 x^3\\
		&+x^2 x^4 (\xi_{11}-2 \xi_{22}+\xi_{33})+3 x^3 (\xi_{13} x^3+\xi_{12} x^4),\\
		(-\frac{3}{2}r^6h^-)_{33}&=-\left((x^1)^2 (\xi_{11}-2 \xi_{22}+\xi_{33})\right)+6 x^1 (\xi_{23} x^2+\xi_{12} x^4)\\
		&-(x^2)^2 (\xi_{11}+\xi_{22}-2 \xi_{33})+6 \xi_{13} x^2 x^4+(x^4)^2 (2 \xi_{11}-\xi_{22}-\xi_{33}),\\
		(-\frac{3}{2}r^6h^-)_{34}&=-3 x^3 (\xi_{12} x^1+\xi_{13} x^2)+3 x^2 (x^1 (\xi_{22}-\xi_{33})+\xi_{12} x^4)\\
		&-3 x^1 (\xi_{23} x^1+\xi_{13} x^4)+3 \xi_{23} (x^2)^2+x^3 x^4 (-2 \xi_{11}+\xi_{22}+\xi_{33}),\\
		(-\frac{3}{2}r^6h^-)_{44}&=-\left((x^1)^2 (\xi_{11}+\xi_{22}-2 \xi_{33})\right)+6 x^1 (\xi_{13} x^3-\xi_{23} x^2)\\
		&-(x^2)^2 (\xi_{11}-2 \xi_{22}+\xi_{33})-6 \xi_{12} x^2 x^3+(x^3)^2 (2 \xi_{11}-\xi_{22}-\xi_{33}).
	\end{align}
	Using the basic metric structure, we use the standard cut-off trick to smoothly extend $r$ from $X\backslash U$ to $X$ such that $r\geq 1$ throughout $X$. We then define the weighted norms.
	\begin{definition}
		\label{defnorm}
		For any $\nu\in\R$, we define the Sobolev norms as follows:
		\begin{align}
			&||\omega||^2_{L^2_\nu(X)}:=\int_{X}|\omega|^2 r^{-4-2\nu}d\Vol_X, &||\omega||^2_{W^{k,2}_\nu(X)}:=\sum_{m=0}^{k}||\nabla^m \omega||^2_{L^2_{\nu-m}(X)},
		\end{align}
		where $\omega$ is a tensor field on $X$.
	\end{definition}
	
	Our proof relies on the following results, which are well-known to experts in the field. For example, see \cite{Mel} for the proof. Note that in a slightly different setting, the first author, Viaclovsky and Zhang have provided a self-contained proof of a similar result, see Proposition 4.5 in \cite{CVZ}.
	
	\begin{proposition}
		\label{wei}
		\label{Weighted-analysis-X}
		Let $(X, g)$ be an ALE manifold of order 4. Then the following properties hold.
		\begin{enumerate}
			
			\item  For any $\nu  \in \mathbb{R}\setminus\mathbb{Z}$ and $k\in\mathbb{N}$, there exist constants $R(X,\nu)>0$ and $C(X, \nu, k)>0$ such that for any $p$-form $\omega\in W^{k+2,2}_\nu(X)$,
			\begin{align}
				\|\omega\|_{W^{k+2,2}_{\nu}(X)} \leq C \cdot (\| \Delta_X \omega\|_{W^{k,2}_{\nu - 2}(X)} + \|\omega\|_{L^2(\{r \le 3R\}\subset X)}).
				\label{e:Weighted-analysis-X}
			\end{align}
			
			\item  For any $\nu  \in \mathbb{R}\setminus\mathbb{Z}$ and $k\in\mathbb{N}$, the operator 
			\begin{equation}
				\Delta_{X}: W^{k+2,2}_\nu (X) \to W^{k,2}_{\nu-2}(X)
			\end{equation}
			is a Fredholm operator.
			Thus for any $p$-form $\omega\in W^{k,2}_{\nu-2}(X)$, 
			\begin{equation}
				\Delta_{X} \tau=\omega
			\end{equation} has a solution $\tau\in W^{k+2,2}_\nu(X)$ if and only if for all $\psi\in \mathcal{H}_{-2-\nu}^p(X)$,
			\begin{align}
				\int_{X} (\omega,\psi)_X d \Vol_X = 0,
			\end{align}
			where $\mathcal{H}_{-2-\nu}^p(X)$ is the space of all harmonic $p$-forms on $X$ in $L^2_{-2-\nu}(X)$. Note that $\mathcal{H}_{-2-\nu}^p(X)\subset W^{k,2}_{-2-\nu}(X)$ for all $k>0$ by standard elliptic estimates.
			
			\item For any $\nu  \in \mathbb{R}\setminus\mathbb{Z}$, $k\in \mathbb{N}$, and $p$-form $\omega\in  W^{k,2}_{\nu}(X)$, there exists some $\tau\in W^{k+2,2}_{\nu+2}(X)$ such that $\Delta_{\R^4} \tau=\omega$ when $r\ge 2R$.
			
			\item  Let $\nu,\nu'  \in \mathbb{R}\setminus\mathbb{Z}$ and $\nu-\nu'\in(0,1)$. Consider any p-form $\omega\in W^{k,2}_{\nu}(\mathbb{R}^4\setminus B_R(0))$ such that $\Delta_{\mathbb{R}^4}\omega=0$ when $r\ge 2R$. If $\mathbb{Z}\cap [\nu',\nu]=\emptyset$, then 
			$\omega \in W^{k,2}_{\nu'}(\mathbb{R}^4\setminus B_R(0))$.
			If there is some $q \in \mathbb{Z}\cap (\nu',\nu)$, then $\omega$ can be written as the sum of a $\mathbb{R}^4$-harmonic form $\sum_{i_1<...<i_p} r^q u_{i_1\cdots i_p}(\theta) dx^{i_1}\wedge...\wedge dx^{i_p}$ and an element in $W^{k,2}_{\nu'}(\mathbb{R}^4\setminus B_R(0))$, where $\theta$ denotes the coordinates on $\mathbb{S}^3$.
		\end{enumerate}
	\end{proposition}
	
	\section{Harmonic functions}
	In this section, we prove Theorem \ref{theoremfunc}. For simplicity, we assume that $\Gamma=\{1\}$, other cases can be reduced to this case by the standard covering space argument.\\
	\subsection{Case 1}
	We first consider the functions $(x^i)^2$, $i=1,2,3,4$. Recall the formula for the Laplacian $\Delta_{X}$. Here, $\Pi^*\Phi^*g=g_{ij}dx^i\otimes dx^j,\ r\geq2R$, and $G=\text{det}(g_{ij})$.
	\begin{align}
		\Delta_X=-\frac{1}{\sqrt{G}}\p_i(\sqrt{G}g^{ij}\p_j),\ \ \Gamma^j_{ji}=\frac{1}{\sqrt{G}}\p_i(\sqrt{G}).
	\end{align}
	We get
	\begin{align}
		\Delta_X((x^i)^2)&=-\left(\left(\Gamma^k_{kj}g^{jl}+\p_j(g^{jl})\right)\p_l+g^{jk}\p_j\p_k\right)(x^i)^2\\
		&=-\left(\Gamma^k_{kj}g^{ji}+\p_j(g^{ji})\right)2x^i-2g^{ii},
	\end{align}
	where we apply the Einstein summation to $j$ and $k$, but not on $i$.\\ 
	In the previous expansion \eqref{expansionmetric}, $h=O(r^{-4}),\ h'=O(r^{-5})$. Note that in this paper, the decaying conditions of the higher derivatives are usually very natural, so we omit them for simplicity, and thus we can perform a Taylor expansion to obtain the leading term:
	\begin{align}
		g^{ij}=\delta_{ij}-h_{ij}-h'_{ij}+O(r^{-8}).
	\end{align}
	Thus
	\begin{align}
		\Gamma^k_{kj}&=\frac{1}{2}g^{kl}(g_{kl|j}+g_{jl|k}-g_{kj|l})\\
		&=\frac{1}{2}(\delta_{kl}-h_{kl}-h'_{kl})\left((h_{kl|j}+h'_{kl|j})+(h_{jl|k}+h'_{jl|k})-(h_{jk|l}+h'_{jk|l})\right)+O(r^{-9})\\
		&=\frac{1}{2}\delta_{kl}\left(h_{kl|j}+h_{jl|k}-h_{jk|l}\right)+O(r^{-6})\\
		&=\frac{1}{2}h_{kk|j}+O(r^{-6}),
	\end{align}
	and
	\begin{align}
		\p_j(g^{ji})=-h_{ji|j}+O(r^{-6}),
	\end{align}
	where we use $T_{I|k}$ to denote $\p_k(T_{I})$ for any tensor $T$, where $I$ is the index.\\
	Thus we get
	\begin{align}
		\Delta_X((x^i)^2)=-\sum_{k=1}^4\left(\frac{1}{2}h_{kk|i}-h_{ki|k}\right)2x^i-2(1-h_{ii})+O(r^{-5}).
	\end{align}
	Note that $\sum_{k=1}^4h_{kk|i}=\sum_{k=1}^4h_{ki|k}=0$. This can be obtained by direct calculation, or by noticing that $h$ is trace-free and divergence-free, see \cite{BH}.\\
	Thus
	\begin{align}
		\Delta_X((x^i)^2)=-2+2h_{ii}+O(r^{-5}),\ r\geq2R.
	\end{align}
	By Proposition \ref{wei}, for any $\nu\notin\Z$, the Laplacian \begin{align}\Delta_X:W^{2,2}_\nu(X)\to L^2_{\nu-2}(X)\end{align}
	is a Fredholm operator. Let $\chi$ be a cutoff function on $X$ such that
	\begin{align}
		\chi= \left\{
		\begin{array}{ll}
			0, & \text{if } r<R, \\
			1, & \text{if } r\geq 2R.
		\end{array}
		\right.
	\end{align}
	We claim that there exists a function $u_{ii}$ on $X$, such that $\chi(x^i)^2-u_{ii} \in W^{2,2}_\nu(X)$ for $\nu$ to be determined, and
	\begin{align}
		\Delta_X(\chi(x^i)^2-u_{ii})=2+\Delta_X(\chi(x^i)^2)=2h_{ii}+O(r^{-5}).
	\end{align}
	That is, $2+\Delta_X(\chi(x^i)^2)\in\image(\Delta_X:W^{2,2}_\nu(X)\to L^2_{\nu-2}(X))$. We note that similar issues for $p$-forms ($p=0,1$) will be encountered, so we address them here.\\
	Firstly it is necessary that
	\begin{align}
		||2+\Delta_X(\chi(x^i)^2)||_{L^2_{\nu-2}(X)}<\infty,
	\end{align}
	which holds if and only if $\nu>-2$.\\
	By Proposition \ref{wei}, it suffices to make $\mathcal{H}^p_{-2-\nu}(X)=0$.
	By standard elliptic estimates, one gets
	\begin{align}
		\mathcal{H}^p_{-2-\nu}(X)=\ker(\Delta_X:W^{k,2}_{-\nu-2}(X)\to W^{k-2,2}_{-\nu-4}(X)),\ \forall k\geq0.
	\end{align}
	In this section, $p=0$. By the maximum principle and elliptic estimates, to ensure that $\ker(\Delta_X)$ trivial, it suffices to make the $L^2-$norm of $\omega$ decay, that is, $\nu>-2$. Thus, in the following, we will take $\nu=-2+\epsilon$, where $0<\epsilon\ll1$. In the next section, $p=1$, we take $\nu=-3+\epsilon$. Then for $H^{p}_{-2-\nu}=H^{p}_{1-\epsilon}$, by Proposition \ref{wei}, the leading order must be $dx^i$. As long as $\Gamma$ is non-trivial, such a thing cannot be $\Gamma-$invariant. So the leading term vanishes, and it decays. Now we apply the Bochner formula
	\begin{align}
		-\frac{1}{2}\Delta_X|\omega|^2&=-\langle\Delta_X\omega,\omega\rangle+|\nabla\omega|^2+Ric(\sharp\omega,\sharp\omega)\\
		&=|\nabla\omega|^2\geq 0,
	\end{align}
	where we used the fact that $\omega$ is harmonic and $X$ is Ricci-flat. This proves the existence of $u_{ii}$ and similar terms.\\
	With the expansion of the metric, one can replace $\Delta_X$ with $\Delta_{\R^4}$, that is, 
	\begin{align}
		\Delta_{\R^4}(\chi(x^i)^2-u_{ii})=2h_{ii}+O(r^{-5}).
	\end{align}
	One can verify that 
	\begin{align}
		\label{defetaii}
		\tilde{\eta}_{ii}=h_{ii}\cdot\frac{r^2}{4}
	\end{align}
	satisfies $\Delta_{\R^4}\tilde{\eta}_{ii}=2h_{ii}$. Therefore 
	\begin{align}\Delta_{\R^4}(\chi(x^i)^2-u_{ii}-\tilde{\eta}_{ii})=O(r^{-5})\in W^{2,2}_{-5+\epsilon},\ r\geq 2R.\end{align}
	By Proposition \ref{wei}, one can find $\tilde{\tilde{\eta}}_{ii}\in W^{2,2}_{-3+\epsilon}$ such that $\Delta_{\R^4}(\tilde{\tilde{\eta}}_{ii})=\Delta_{\R^4}(\chi(x^i)^2-u_{ii}-\tilde{\eta}_{ii})$ for $r\geq 2R.$ That is,
	\begin{align}
		\Delta_{\R^4}((x^i)^2-u_{ii}-\tilde{\eta}_{ii}-\tilde{\tilde{\eta}}_{ii})=0,\ r\geq 2R.
	\end{align}
	By Proposition \ref{wei}, one can obtain
	\begin{align}
		(x^i)^2-u_{ii}-\tilde{\eta}_{ii}-\tilde{\tilde{\eta}}_{ii}&=(r^{-2}\text{ ordered homogeneous harmonic function})+O(r^{-3})\\
		&=-\frac{C_{ii}}{r^2}+O(r^{-3}),\ r\geq 2R.
	\end{align}
	Therefore one finally gets
	\begin{align}
		u_{ii}=(x^i)^2-\tilde{\eta}_{ii}+\frac{C_{ii}}{r^2}+O(r^{-3+\epsilon}),\ r\geq 2R,
	\end{align}
	where we have omitted  $\tilde{\tilde{\eta}}_{ii}$.
	\subsection{Case 2}
	As for the harmonic functions of type $x^ix^j$, the calculation is similar.
	\begin{align}
		\Delta_X(\chi x^ix^j)&=-\left(\Gamma^k_{kl}g^{li}+\p_l(g^{li})\right)x^j-\left(\Gamma^k_{kl}g^{lj}+\p_l(g^{lj})\right)x^i-2g^{ij}+O(r^{-5})\\
		&=\sum_{k=1}^4(\frac{1}{2}h_{kk|i}+h_{ki|k})x^j+\sum_{k=1}^4(\frac{1}{2}h_{kk|j}+h_{kj|k})x^i+2h_{ij}+O(r^{-5})\\
		&=2h_{ij}+O(r^{-5}),\ r\geq 2R.
	\end{align}
	By similar arguments, one can find a function $v_{ij}\in W^{2,2}_{-2+\epsilon}(X)$, such that $\Delta_X(\chi x^ix^j)=\Delta_X(v_{ij}).$ Moreover, one obtains the harmonic function
	\begin{align}
		u_{ij}&=\chi x^ix^j-v_{ij}\\
		&=x^ix^j-\tilde{\eta}_{ij}+\frac{C_{ij}}{r^2}+O(r^{-3+\epsilon}),\ r\geq2R,
	\end{align}
	where 
	\begin{align}
		\label{defetaij}
		\tilde{\eta}_{ij}=h_{ij}\cdot\frac{r^2}{4}
	\end{align}
	satisfies $\Delta_{\R^4}\tilde{\eta}_{ij}=2h_{ij}.$
	\subsection{Determine the constant}
	Let $F$ be $u_{ii}$ or $u_{ij}$ ($i\neq j$). By integrating $\Delta_XF$ on the region $\Omega_{\rho}$ such that $\Pi^*\Phi^*(\Omega_{\rho}-U)=[-\rho,\rho]^4-B_R(0)$, we get
	\begin{align}
		&\int_{\Omega_{\rho}} \Delta_XF d\Vol_X\\
		&=-\int_{\Omega_{\rho}} d*_XdF\\
		&=-\int_{\p\Omega_{\rho}} *_XdF.
	\end{align}
	With the convention that
	\begin{align}
		\omega\wedge*_X\eta=\langle\omega,\eta\rangle d\Vol_X,\ \ \langle dx^i,dx^j\rangle=g^{ij},
	\end{align}
	we get
	\begin{align}
		*_XdF=\sum_{j,k=1}^4\frac{\p F}{\p x^k}g^{kj}\sqrt{G}(-1)^{j-1}dx^1\wedge\cdots\wedge\widehat{dx^j}\wedge\cdots\wedge dx^4,\ r\geq 2R,
	\end{align}
	where $G=\det(g_{ij})=1+\sum_i h_{ii}+O(r^{-5})=1+O(r^{-5})$ because $h$ is trace-free. Combine this with $g^{ij}=\delta_{ij}-h_{ij}+O(r^{-5})$, one gets
	\begin{align}
		\int_{\p\Omega_{\rho}} *_XdF &= \sum_{j=1}^{4}\int_{\p[-\rho,\rho]^4}\left(\frac{\p F}{\p x^j}-\sum_{k=1}^{4}h_{kj}\frac{\p F}{\p x^k}\right)(-1)^{j-1}dx^1\wedge\cdots\wedge\widehat{dx^j}\wedge\cdots\wedge dx^4\\
		&+O(\rho^{-1+\epsilon}).
	\end{align}
	\subsection{Case 1}
	We choose $F=u_{ii}=\chi(x^i)^2-\tilde{\eta}_{ii}+\frac{C_{ii}}{r^2}+O(r^{-3+\epsilon})$. On the LHS, we get
	\begin{align}
		\int_{\Omega_{\rho}}\Delta_X(u_{ii})d\Vol_X=-2\Vol_X(\Omega_{\rho}).
	\end{align}
	Denote $\gamma_j=(-1)^{j-1}dx^1\wedge\cdots\wedge\widehat{dx^j}\wedge\cdots\wedge dx^4$. On the RHS, we get
	\begin{align}
		&\int_{\Omega_{\rho}}\left(\Delta_X(\chi(x^i)^2-\tilde{\eta}_{ii}+\frac{C_{ii}}{r^2})+O(r^{-5+\epsilon})\right)d\Vol_X\\
		&=-\int_{\p \Omega_{\rho}}*_Xd((x^i)^2-\tilde{\eta}_{ii}+\frac{C_{ii}}{r^2})+O(\rho^{-1+\epsilon})\\
		&=-\sum_{j=1}^4\int_{\p[-\rho,\rho]^4}\left(\frac{\p}{\p x^j}((x^i)^2-\tilde{\eta}_{ii}+\frac{C_{ii}}{r^2})-\sum_{k=1}^4h_{kj}\frac{\p}{\p x^j}((x^i)^2)\right)\gamma_j+O(\rho^{-1+\epsilon})\\
		&=4\pi^2C_{ii}-2\Vol_{g_0}([-\rho,\rho]^4)+\sum_{j=1}^4\int_{\p[-\rho,\rho]^4}\left(\frac{\p}{\p x^j}(\tilde{\eta}_{ii})+\sum_{k=1}^4h_{kj}\frac{\p}{\p x^j}((x^i)^2)\right)\gamma_j\\
		&+O(\rho^{-1+\epsilon}),
	\end{align}
	where we have used
	\begin{align}
		\sum_{j=1}^4\int_{\p[-\rho,\rho]^4}\frac{\p}{\p x^j}(\frac{1}{r^2})\gamma_j=-4\pi^2,\ &\sum_{j=1}^4\int_{\p[-\rho,\rho]^4}\frac{\p}{\p x^j}((x^i)^2)\gamma_j=2\Vol_{g_0}([-\rho,\rho]^4).
	\end{align}
	Notice that the leading terms of the integrand $\frac{\p}{\p x^j}((x^i)^2-\tilde{\eta}_{ii}+\frac{C_{ii}}{r^2})-\sum_{k=1}^4h_{kj}\frac{\p}{\p x^j}((x^i)^2)$ are all of order $r^{-3}$. Thus, letting $\rho\to\infty$, the errors vanish, and the remaining integral is invariant. We may thus replace $\Omega_{\rho}$ with $\Omega=\Omega_1.$\\
	Now we compute,
	\begin{align}
		&\sum_{j=1}^4\int_{\p[-\rho,\rho]^4}\left(\frac{\p}{\p x^j}(\tilde{\eta}_{ii})+\sum_{k=1}^4h_{kj}\frac{\p}{\p x^j}((x^i)^2)\right)\gamma_j\\
		&=\sum_{j=1}^4\sum_{l=1}^4\int_{\Sigma_l^+ +\Sigma_l^-}\left(\frac{\p}{\p x^j}(\tilde{\eta}_{ii})+\sum_{k=1}^4h_{kj}\frac{\p}{\p x^j}((x^i)^2)\right)\gamma_j\\
		&=\sum_{j=1}^4\int_{\Sigma_j^+ +\Sigma_j^-}\left(\frac{\p}{\p x^j}(\tilde{\eta}_{ii})+\sum_{k=1}^4h_{kj}\frac{\p}{\p x^j}((x^i)^2)\right)\gamma_j\\
		&=\sum_{j=1}^4\int_{\Sigma_j^+}\left(\frac{\p}{\p x^j}(\tilde{\eta}_{ii})+\sum_{k=1}^4h_{kj}\frac{\p}{\p x^j}((x^i)^2)\right)dx^1\cdots\widehat{dx^j}\cdots dx^4\\
		&-\sum_{j=1}^4\int_{\Sigma_j^-}\left(\frac{\p}{\p x^j}(\tilde{\eta}_{ii})+\sum_{k=1}^4h_{kj}\frac{\p}{\p x^j}((x^i)^2)\right)dx^1\cdots\widehat{dx^j}\cdots dx^4,
	\end{align}
	where $\Sigma_l^\pm=\{x^l=\pm1\}\times[-1,1]^3$, e.g. $\Sigma_1^+=\{x^1=1\}\times[-1,1]^3$.\\
	Then we can find that the integrals equal zero by direct calculation, and the Wolfram Mathematica program is attached in the appendix.\\
	Thus we finally get
	\begin{align}
		4\pi^2C_{ii}=\lim_{\rho\to \infty}2(\Vol_{g_0}([-\rho,\rho]^4)-\Vol_X(\Omega_{\rho}))=:-2\mathcal{V},
	\end{align}
	where $\mathcal{V}$ is the renormalized volume, which is finite according to Biquard and Hein's works \cite{BH}. Note that when $\Gamma\neq\{1\}$, there will be a factor $|\Gamma|$ for the covering reason.
	\subsection{Case 2}Now we choose $F=u_{ij}=x^ix^j-\tilde{\eta}_{ij}+\frac{C_{ij}}{r^2}+O(r^{-3+\epsilon})$.\\
	Notice that the integral of the 1-order term $\sum_{j=1}^{4}\frac{\p }{\p x^k}(x^ix^j)\gamma_j$ is obviously zero, and all the other leading terms are of -3-order, thus by a similar procedure, we get
	\begin{align}
		2\pi^2 C_{ij}&=\sum_{j=1}^4\int_{\Sigma_j^+}\left(\frac{\p}{\p x^j}(\tilde{\eta}_{ij})+\sum_{k=1}^{4}h_{kj}\frac{\p}{\p x^k}(x^ix^j)\right)dx^1\cdots\widehat{dx^j}\cdots dx^4\\
		&-\sum_{j=1}^4\int_{\Sigma_j^-}\left(\frac{\p}{\p x^j}(\tilde{\eta}_{ij})+\sum_{k=1}^{4}h_{kj}\frac{\p}{\p x^k}(x^ix^j)\right)dx^1\cdots\widehat{dx^j}\cdots dx^4.
	\end{align}
	Thus one can determine $C_{ij}$ by calculating the integral on the RHS. In fact, our Wolfram Mathematica program tells us that 
	\begin{align}
		C_{ij}=0\ (i\neq j).
	\end{align}
	\section{Harmonic 1-forms}
	In this part, we consider the (invariant) harmonic 1-forms on ALE manifolds. Specifically, we only need to consider $\omega_{i_1i_2}$, which is the 1-form asymptotic to $x^{i_1}dx^{i_2}$. For simplicity, we still assume that $\Gamma=\{1\}$.
	\subsection{Find the expansions}
	The process is similar to that described in section 2.1, but it requires more calculation. Before proceeding, let's recall our convention for the Hodge stars, i.e.
	\begin{align}
		*_X(dx^{i_1}\wedge\cdots\wedge dx^{i_k})=\sum_{\substack{j_{k+1}<\cdots<j_n,\\j_r\neq j_s\ for\ r\neq s}}\sqrt{G}g^{i_1j_1}\cdots g^{i_kj_k}\varepsilon_{j_1\cdots j_n}dx^{j_{k+1}}\wedge\cdots\wedge dx^{j_n},
	\end{align}
	where $\varepsilon_{j_1\cdots j_n}$ is the Levi-Civita symbol with $\varepsilon_{1\cdots n}=1$ and $\varepsilon_{j_1\cdots j_n}\neq0$ if and only if all $j_r$ are different with each other. We begin with Proposition \ref{wei} to get $\tilde{\omega}_{i_1i_2}\in W^{2,2}_{-3+\epsilon}(X)$ that satisfies
	\begin{align}
		\Delta_X\tilde{\omega}_{i_1i_2}=\Delta_X(\chi x^{i_1}dx^{i_2})=(d\delta_X+\delta_X d)(\chi x^{i_1}dx^{i_2}).
	\end{align}
	First, we note that, for a 1-form $\tilde{\omega}=\tilde{\omega}_idx^i$, $\delta_X\tilde{\omega}=-\frac{1}{\sqrt{G}}\frac{\p}{\p x^j}(\tilde{\omega}_ig^{ij}\sqrt{G})$. Thus
	\begin{align}
		\delta_X (\chi x^{i_1}dx^{i_2})&=-\frac{1}{\sqrt{G}}\frac{\p}{\p x^j}(\chi x^{i_1}g^{i_2j}\sqrt{G})\\
		&=-\frac{\p}{\p x^j}(\chi x^{i_1}g^{i_2j})+O(r^{-5})\\
		&=-\delta_{i_1i_2}+h_{i_1i_2}+\sum_{j=1}^4x^{i_1}h_{i_2j|j}+O(r^{-5}),\ r\geq 2R,
	\end{align}
	and therefore
	\begin{align}
		d\delta_X (\chi x^{i_1}dx^{i_2})&=\sum_{k=1}^4\frac{\p}{\p x^k}(h_{i_1i_2}+\sum_{j=1}^4x^{i_1}h_{i_2j|j})dx^k+O(r^{-6}),\ r\geq 2R.
	\end{align}
	For the second term, we note that for a 2-form $\tilde{\omega}=\sum_{i,j}\tilde{\omega}_{ij}dx^i\wedge dx^j$, we have\\ $\delta_X\tilde{\omega}=\frac{1}{\sqrt{G}}\frac{\p}{\p x^j}\left((\tilde{\omega}_{kl}-\tilde{\omega}_{lk})g^{ik}g^{jl}\sqrt{G}\right)g_{im}dx^{m}$. Thus
	\begin{align}
		\delta_X d (x^{i_1}dx^{i_2})&=\frac{1}{\sqrt{G}}\frac{\p}{\p x^j}\left((\delta_{i_1k}\delta_{i_2l}-\delta_{i_1l}\delta_{i_2k})g^{ik}g^{jl}\sqrt{G}\right)g_{im}dx^{m}\\
		&=-\sum_{i,j,k,l=1}^4(\delta_{i_1k}\delta_{i_2l}-\delta_{i_1l}\delta_{i_2k})\frac{\p}{\p x^j}\left(\delta_{ik}h_{jl}+h_{ik}\delta_{jl}\right)dx^i+O(r^{-6}),\ r\geq 2R.
	\end{align}
	As a result, we get
	\begin{align}
		\Delta_X(\chi x^{i_1}dx^{i_2})&=\sum_{k=1}^4\frac{\p}{\p x^k}(h_{i_1i_2}+\sum_{j=1}^4x^{i_1}h_{i_2j|j})dx^k\\
		&+\sum_{i,j,k,l=1}^4(\delta_{i_1l}\delta_{i_2k}-\delta_{i_1k}\delta_{i_2l})\frac{\p}{\p x^j}\left(\delta_{ik}h_{jl}+h_{ik}\delta_{jl}\right)dx^i+O(r^{-6})\\
		&:=L_{i_1i_2}+O(r^{-6}),\ r\geq 2R.
	\end{align}
	One can verify that
	\begin{align}
		\label{defmuij}
		\tilde{\mu}_{i_1i_2}=L_{i_1i_2}\cdot\frac{r^2}{12}
	\end{align}
	satisfies $\Delta_{\R^4}\tilde{\mu}_{i_1i_2}=L_{i_1i_2}$.\\
	By Proposition \ref{wei}, we get the expansion of harmonic 1-form
	\begin{align}
		\omega_{i_1i_2}=\chi x^{i_1}dx^{i_2}-\tilde{\omega}_{i_1i_2},\ r\geq 2R,
	\end{align}
	where
	\begin{align}
		\tilde{\omega}_{i_1i_2}:=\tilde{\mu}_{i_1i_2}+\sum_{k,l=1}^4\frac{\Con[i_1,i_2,k,l] x^k}{r^4}dx^{l}
	\end{align}
	\subsection{Determine the relation of the constants} In this section, we will derive the restriction equations using these methods:
	\begin{enumerate}
		\item Differentiation of harmonic functions.
		\item Divergence arguments.
		\item Integral of the Laplacian of forms.
		\item Integral of the covariant derivative of forms.
	\end{enumerate}
	First, recall our expansion:
	\begin{align}
		u_{ii}&=(x^i)^2-\tilde{\eta}_{ii}+\frac{C_{ii}}{r^2}+O(r^{-3+\epsilon}),\ r\geq 2R,\\
		C_{ii}&=-\frac{\mathcal{V}}{2\pi^2},\quad(i=1,2,3,4)\\
		\omega_{i_1i_2}&=\chi x^{i_1}dx^{i_2}-\tilde{\omega}_{i_1i_2}\\
		&=x^{i_1}dx^{i_2}-\tilde{\mu}_{i_1i_2}-\sum_{k,l=1}^4\frac{\Con[i_1,i_2,k,l] x^k}{r^4}dx^{l}+O(r^{-4+\epsilon}),\ r\geq 2R.
	\end{align}
	We differentiate the first equation and obtain
	\begin{align}
		du_{ii}&=d\left((x^i)^2-\tilde{\eta}_{ii}+\frac{C_{ii}}{r^2}\right)+O(r^{-4+\epsilon}),\ r\geq 2R.
	\end{align}
	Thus, $du_{ii}$ is a harmonic 1-form asymptotic to $2x^idx^i$. We see that $du_{ii}-2\omega_{ii}$ is a decaying harmonic 1-form, which vanishes by the Ricci-flat condition. Thus, we obtain
	\begin{align}
		du_{ii}=2\omega_{ii},
	\end{align}
	that is,
	\begin{align}
		2\left(\tilde{\mu}_{ii}+\sum_{k,l=1}^4\frac{\Con[i,i,k,l] x^k}{r^4}dx^{l}\right)+d\left(-\tilde{\eta}_{ii}+\frac{C_{ii}}{r^2}\right)=0,\ r\geq 2R.
	\end{align}
	Similarly, for $i\neq j$, we have
	\begin{align}
		du_{ij}=\omega_{ij}+\omega_{ji}.
	\end{align}
	Next, we consider divergence arguments. For the harmonic 1-form $\omega_{i_1i_2}=\chi x^{i_1}dx^{i_2}-\tilde{\omega}_{i_1i_2}$, $\delta_X\omega_{i_1i_2}$ is a decaying  harmonic function plus $-\delta_{i_1i_2}$. By applying the maximum principle, $\delta_X\omega_{i_1i_2}=-\delta_{i_1i_2}$. Since $\tilde{\omega}_{i_1i_2}$ is of order -3, one can calculate the divergence of $\tilde{\omega}_{i_1i_2}$ using the Euclidean divergence. That is,
	\begin{align}
		\delta_X\omega_{i_1i_2}+\delta_{i_1i_2}&=\delta_X(\chi x^{i_1}dx^{i_2}-\tilde{\omega}_{i_1i_2})+\delta_{i_1i_2}\ (r\geq 2R)\\
		&=\sum_{l=1}^4\frac{\p}{\p x^l}(\tilde{\omega}_{i_1i_2;l})+h_{i_1i_2}+\sum_{j=1}^4x^{i_1}h_{i_2j|j}+O(r^{-5+\epsilon})\\
		&=0.
	\end{align}
	Next, we consider integrals involving the Laplacian. For two harmonic 1-forms $\omega_{i_1i_2}$ and $\omega_{i_3i_4}$,using the divergence arguments, we have
	\begin{align}
		0&=\int_{\Omega_{\rho}}\langle\Delta_X\omega_{i_1i_2},\omega_{i_3i_4}\rangle d\Vol_X-\int_{\Omega_{\rho}}\langle\Delta_X\omega_{i_3i_4},\omega_{i_1i_2}\rangle d\Vol_X\\
		&=(\int_{\p\Omega_{\rho}}*_Xd\omega_{i_1i_2}\wedge\omega_{i_3i_4}+\int_{\p\Omega_{\rho}}*_Xd*_X\omega_{i_1i_2}\wedge*_X\omega_{i_3i_4})\\
		&-(\int_{\p\Omega_{\rho}}*_Xd\omega_{i_3i_4}\wedge\omega_{i_1i_2}+\int_{\p\Omega_{\rho}}*_Xd*_X\omega_{i_3i_4}\wedge*_X\omega_{i_1i_2})\\
		&=\int_{\p\Omega_{\rho}}*_Xd\omega_{i_1i_2}\wedge\omega_{i_3i_4}-\int_{\p\Omega_{\rho}}*_Xd\omega_{i_3i_4}\wedge\omega_{i_1i_2}\\
		&=\int_{\p\Omega_{\rho}}\alpha_{i_1i_2i_3i_4}-\int_{\p\Omega_{\rho}}\alpha_{i_3i_4i_1i_2}.
	\end{align}
	Here we have used the leading term arguments, and the $\alpha_{i_1i_2i_3i_4}$ denotes the leading term of $*_Xd\omega_{i_1i_2}\wedge\omega_{i_3i_4}$. Furthermore, one can calculate $\alpha_{i_1i_2i_3i_4}$ explicitly, say,
	\begin{align}
		&\alpha_{i_1i_2i_3i_4}\\
		&=-\sum_{\substack{j_3<j_4,\\j_r\neq j_s\ for\ r\neq s}}\sum_{t=1}^4(\frac{\p}{\p x^{j_1}}(\tilde{\omega}_{i_1i_2;j_2})x^{i_3}\delta_{i_4 t}+\delta_{i_1j_1}\delta_{i_2j_2}\tilde{\omega}_{i_3i_4;t}\\
		&+h_{i_1j_1}\delta_{i_2j_2}x^{i_3}\delta_{i_4 t}+\delta_{i_1j_1}h_{i_2j_2}x^{i_3}\delta_{i_4 t})\varepsilon_{j_1j_2j_3j_4}dx^{j_3}\wedge dx^{j_4}\wedge dx^t.
	\end{align}
	Finally, we consider the integral of covariant derivatives. Let $\omega_{i_1i_2}$ and $\omega_{i_3i_4}$ be two harmonic 1-forms.\\
	Recall that we use the following convention:
	\begin{align}
		\langle dx^{i_1}\otimes\cdots\otimes dx^{i_k},dx^{j_1}\otimes\cdots\otimes dx^{j_k}\rangle=\frac{1}{k!}g^{i_1j_1}\cdots g^{i_kj_k}.
	\end{align}
	Furthermore, for $\phi=\phi_i dx^i$ and $Y=Y^i\p_i$, we have:
	\begin{align}
		\nabla\phi&=\phi_{i,j}dx^j\otimes dx^i=(\p_j\phi_i-\phi_k\Gamma^k_{ij})dx^j\otimes dx^i,\\
		\nabla Y&=Y^i_{\ ,j}dx^j\otimes\p_i=(\p_j Y^i+Y^k\Gamma^i_{kj})dx^j\otimes \p_i.
	\end{align}
	Now we consider the following integral, which we aim to express as a boundary integral:
	\begin{align}
		&\int_{\Omega_{\rho}}\langle\nabla\omega_{i_1i_2},\nabla\omega_{i_3i_4}\rangle d\Vol_X.
	\end{align}
	Define the 1-form $\beta_{i_1 i_2 i_3 i_4}$ as
	\begin{align}
		\beta_{i_1 i_2 i_3 i_4}(Y):=\frac{1}{2}\langle\omega_{i_1i_2},\nabla_{Y}\omega_{i_3i_4}\rangle.
	\end{align}
	Then we claim that
	\begin{align}
	&\int_{\Omega_{\rho}}\langle\nabla\omega_{i_1i_2},\nabla\omega_{i_3i_4}\rangle d\Vol_X\\
	&=\int_{\Omega_{\rho}} \text{div}(\beta_{i_1 i_2 i_3 i_4}) d\Vol_X\\
	&=\int_{\p\Omega_{\rho}}\beta_{i_1 i_2 i_3 i_4}(N)(N\lrcorner d\Vol_X),
	\end{align}
	where $N$ denotes the exterior unit normal vector field of $\Omega_{\rho}$.\\
	To see this, we compute the divergence in local coordinates:
	\begin{align}
		\text{div}(\beta_{i_1 i_2 i_3 i_4})&=g^{ij}(\nabla\beta_{i_1 i_2 i_3 i_4})_{ij},\\
		(\beta_{i_1 i_2 i_3 i_4})_i&=\frac{1}{2}\langle(\omega_{i_1i_2})_k dx^k,\nabla_{i}((\omega_{i_3i_4})_l dx^l)\rangle\\
		&=\frac{1}{2}(\omega_{i_1i_2})_k(\omega_{i_3i_4})_{l,i}g^{kl}, \ r\geq 2R,\\
		g^{ij}(\nabla\beta_{i_1 i_2 i_3 i_4})_{ij}&=g^{ij}(\beta_{i_1 i_2 i_3 i_4})_{i,j}\\
		&=\frac{1}{2}g^{ij} g^{kl} ((\omega_{i_1i_2})_{k,j}(\omega_{i_3i_4})_{l,i}+(\omega_{i_1i_2})_k(\omega_{i_3i_4})_{l,ij}),\ r\geq 2R.
	\end{align}
	In the last equation we have used the fact that metric tensor is parallel. The first summation is exactly $\langle\nabla\omega_{i_1i_2},\nabla\omega_{i_3i_4}\rangle$. On the other hand, from Weitzenb\"ock formula and that $X$ is Ricci-flat, we found the second term is $\frac{1}{2}\langle\omega_{i_1i_2},-\Delta_X\omega_{i_3i_4}\rangle$, thus vanishes because $\omega_{i_3i_4}$ is harmonic.\\
	Now we can do the calculation in local coordinates. For $\Sigma_j^\pm$, we denote their exterior unit normal vector fields by $N_j^{\pm}=\pm g^{ij}(g^{jj})^{-\frac{1}{2}} \p_i$ for $j=1,\cdots,4$. Then
	\begin{align}
		N^\pm_j\lrcorner d\Vol_X&=\pm\sqrt{G}\varepsilon_{i_1 i_2 i_3 i_4}
		 dx^{i_1}(g^{ij}(g^{jj})^{-\frac{1}{2}}\p_i) dx^{i_2}\otimes dx^{i_3}\otimes dx^{i_4}\\
		 &=\pm(-1)^{i-1}g^{ij}(g^{jj})^{-\frac{1}{2}}dx^{1}\wedge\cdots\wedge\widehat{dx^i}\wedge\cdots\wedge dx^4+O(r^{-5})\\
		 &=\pm g^{ij}(g^{jj})^{-\frac{1}{2}}\gamma_i+O(r^{-5}).
	\end{align}
	On the other hand, $\flat N^\pm_j=\pm (g^{jj})^{-\frac{1}{2}}dx^j$ for $j=1,\cdots,4$. Therefore,
	\begin{align}
		\beta_{i_1 i_2 i_3 i_4}(N^\pm_j)&=\langle\beta_{i_1 i_2 i_3 i_4},\flat N^\pm_j\rangle\\
		&=\pm\frac{1}{2}(\omega_{i_1i_2})_k(\omega_{i_3i_4})_{l,i}g^{kl}g^{ij}(g^{jj})^{-\frac{1}{2}}.
	\end{align}
	Thus, using the leading term argument to eliminate the lower order terms as before, we will get
	\begin{align}
		&\int_{\Omega_{\rho}}\langle\nabla\omega_{i_1i_2},\nabla\omega_{i_3i_4}\rangle d\Vol_X\\
		&=\int_{\p\Omega_{\rho}}\beta_{i_1 i_2 i_3 i_4}(N)(N\lrcorner d\Vol_X)\\
		&=\sum_{j=1}^4 \int_{\Sigma_j^+}\beta_{i_1 i_2 i_3 i_4}(N^+_j)(N^+_j\lrcorner d\Vol_X)+\sum_{j=1}^4 \int_{\Sigma_j^-}\beta_{i_1 i_2 i_3 i_4}(N^-_j)(N^-_j\lrcorner d\Vol_X)\\
		&=\sum_{j=1}^4 \int_{\Sigma_j^+}\beta_{i_1 i_2 i_3 i_4}(N^+_j)(N^+_j\lrcorner d\Vol_X)+\sum_{j=1}^4 \int_{\Sigma_j^-}\beta_{i_1 i_2 i_3 i_4}(N^+_j)(N^+_j\lrcorner d\Vol_X)\\
		&=\sum_{j=1}^{4} \int_{\Sigma_j^+}\frac{1}{2}(\omega_{i_1i_2})_k(\omega_{i_3i_4})_{l,i}g^{kl}g^{ij}g^{mj}(g^{jj})^{-1}\gamma_m\\
		&+\int_{\Sigma_j^-}\frac{1}{2}(\omega_{i_1i_2})_k(\omega_{i_3i_4})_{l,i}g^{kl}g^{ij}g^{mj}(g^{jj})^{-1}\gamma_m\\
		&=\sum_{j=1}^{4}\int_{[-1,1]^3}\frac{1}{2}(\omega_{i_1i_2})_k(\omega_{i_3i_4})_{l,i}g^{kl}g^{ij}|^{x^j=1}_{x^j=-1}dx^1\cdots \widehat{dx^j}\cdots dx^4\\
		&=\sum_{j=1}^{4}\int_{[-1,1]^3}-\frac{1}{2}[(\tilde{\omega}_{i_1i_2})_{i_4}\delta_{i_3j}+x^{i_1}(\p_j (\tilde{\omega}_{i_3i_4})_{i_2}+x^{i_3}\Gamma^{i_4}_{i_2j}\\
		&\left.+\delta_{i_2i_4}h_{i_3j}+\delta_{i_3j}h_{i_2i_4})]\right|^{x^j=1}_{x^j=-1}dx^1\cdots \widehat{dx^j}\cdots dx^4.
	\end{align}
	
	By swapping the indexes, we get
	\begin{align}
		\int_{\p\Omega_{\rho}}\beta_{i_1 i_2 i_3 i_4}(N)(N\lrcorner d\Vol_X)=\int_{\p\Omega_{\rho}}\beta_{i_3 i_4 i_1 i_2}(N)(N\lrcorner d\Vol_X).
	\end{align}
	All of these calculation can be done by Mathematica.\\
	Then we will gain an enormous collection of linear equations. Solve it and we will get some restrictions about $\zeta_{ij}$, $\xi_{ij}$ and $\Con[i,j,k,l]$.
	\section{The Constants}
	The variables can be categorized into 7 types by symmetry:
	\begin{align}
		&\text{I}.\Con[i,i,i,i];\\
		&\text{II}.\Con[i,i,i,j],\Con[i,i,j,i],\Con[i,j,i,i],\Con[j,i,i,i];\\
		&\text{III}.\Con[i,i,j,j];\\
		&\text{IV}.\Con[i,j,i,j],\Con[i,j,j,i];\\
		&\text{V}.\Con[i,i,j,k],\Con[j,k,i,i];\\
		&\text{VI}.\Con[i,j,k,i],\Con[j,i,i,k],\Con[i,j,i,k],\Con[j,i,k,i];\\
		&\text{VII}.\Con[i,j,k,l].
	\end{align}
\subsection{Type I}
\begin{align}
	\Con[i,i,i,i] = -\frac{\mathcal{V}}{2\pi^2}\ (i=1,2,3,4);
\end{align}
\subsection{Type II}
\begin{align}
	\Con[i,i,i,j] &= \Con[i,i,j,i] = \Con[i,j,i,i] = \Con[j,i,i,i] = 0\ (i,j=1,2,3,4,\ i\neq j); 
\end{align}
\subsection{Type III}
\begin{align}
	\Con[1, 1, 2, 2] &= \Con[2, 2, 1, 1] = \Con[3, 3, 4, 4] = \Con[4, 4, 3, 3] \nonumber \\
	&= -\frac{\mathcal{V}}{2\pi^2} - \frac{\zeta_{11}}{9} + \frac{\zeta_{22}}{18} + \frac{\zeta_{33}}{18} - \frac{\xi_{11}}{9} + \frac{\xi_{22}}{18} + \frac{\xi_{33}}{18}; \\
	\Con[1, 1, 3, 3] &= \Con[2, 2, 4, 4] = \Con[3, 3, 1, 1] = \Con[4, 4, 2, 2] \nonumber \\
	&= -\frac{\mathcal{V}}{2\pi^2} + \frac{\zeta_{11}}{18} - \frac{\zeta_{22}}{9} + \frac{\zeta_{33}}{18} + \frac{\xi_{11}}{18} - \frac{\xi_{22}}{9} + \frac{\xi_{33}}{18}; \\
	\Con[1, 1, 4, 4] &= \Con[2, 2, 3, 3] = \Con[3, 3, 2, 2] = \Con[4, 4, 1, 1] \nonumber \\
	&= -\frac{\mathcal{V}}{2\pi^2} + \frac{\zeta_{11}}{18} + \frac{\zeta_{22}}{18} - \frac{\zeta_{33}}{9} + \frac{\xi_{11}}{18} + \frac{\xi_{22}}{18} - \frac{\xi_{33}}{9};
\end{align}
\subsection{Type IV}
\begin{align}
	\Con[i, j, i, j] &= \Con[j, i, j, i];\ (i,j=1,2,3,4,\ i\neq j)\\
	\Con[1, 2, 2, 1] &= \Con[2, 1, 1, 2] \nonumber \\
	&= \frac{\zeta_{11}}{9} - \frac{\zeta_{22}}{18} - \frac{\zeta_{33}}{18} + \frac{\xi_{11}}{9} - \frac{\xi_{22}}{18} - \frac{\xi_{33}}{18} - \Con[1, 2, 1, 2]; \\
	\Con[3, 4, 4, 3] &= \Con[4, 3, 3, 4] \nonumber \\
	&= \frac{\zeta_{11}}{9} - \frac{\zeta_{22}}{18} - \frac{\zeta_{33}}{18} + \frac{\xi_{11}}{9} - \frac{\xi_{22}}{18} - \frac{\xi_{33}}{18} - \Con[3, 4, 3, 4]; \\
	\Con[1, 3, 3, 1] &= \Con[3, 1, 1, 3] \nonumber \\
	&= -\frac{\zeta_{11}}{18} + \frac{\zeta_{22}}{9} - \frac{\zeta_{33}}{18} - \frac{\xi_{11}}{18} + \frac{\xi_{22}}{9} - \frac{\xi_{33}}{18} - \Con[1, 3, 1, 3]; \\
	\Con[1, 4, 4, 1] &= \Con[4, 1, 1, 4] \nonumber \\
	&= -\frac{\zeta_{11}}{18} - \frac{\zeta_{22}}{18} + \frac{\zeta_{33}}{9} - \frac{\xi_{11}}{18} - \frac{\xi_{22}}{18} + \frac{\xi_{33}}{9} - \Con[1, 4, 1, 4]; \\
	\Con[2, 3, 3, 2] &= \Con[3, 2, 2, 3] \nonumber \\
	&= -\frac{\zeta_{11}}{18} - \frac{\zeta_{22}}{18} + \frac{\zeta_{33}}{9} - \frac{\xi_{11}}{18} - \frac{\xi_{22}}{18} + \frac{\xi_{33}}{9} - \Con[2, 3, 2, 3]; \\
	\Con[2, 4, 4, 2] &= \Con[4, 2, 2, 4] \nonumber \\
	&= -\frac{\zeta_{11}}{18} + \frac{\zeta_{22}}{9} - \frac{\zeta_{33}}{18} - \frac{\xi_{11}}{18} + \frac{\xi_{22}}{9} - \frac{\xi_{33}}{18} - \Con[2, 4, 2, 4];
\end{align}
\subsection{Type V}
\begin{align}
	\Con[1, 1, 2, 3] &= \Con[1, 1, 3, 2] = \Con[2, 3, 1, 1] = \Con[3, 2, 1, 1] = -\frac{\zeta_{12}}{18} - \frac{\xi_{12}}{18}; \\
	\Con[1, 1, 2, 4] &= \Con[1, 1, 4, 2] = \Con[2, 4, 1, 1] = \Con[4, 2, 1, 1] = -\frac{\zeta_{13}}{18} - \frac{\xi_{13}}{18}; \\
	\Con[1, 1, 3, 4] &= \Con[1, 1, 4, 3] = \Con[3, 4, 1, 1] = \Con[4, 3, 1, 1] = -\frac{\zeta_{23}}{18} - \frac{\xi_{23}}{18}; \\
	\Con[1, 2, 3, 3] &= \Con[2, 1, 3, 3] = \Con[3, 3, 1, 2] = \Con[3, 3, 2, 1] = -\frac{\zeta_{23}}{18} + \frac{\xi_{23}}{18}; \\
	\Con[1, 2, 4, 4] &= \Con[2, 1, 4, 4] = \Con[4, 4, 1, 2] = \Con[4, 4, 2, 1] = \frac{\zeta_{23}}{18} - \frac{\xi_{23}}{18}; \\
	\Con[1, 3, 2, 2] &= \Con[2, 2, 1, 3] = \Con[2, 2, 3, 1] = \Con[3, 1, 2, 2] = \frac{\zeta_{13}}{18} - \frac{\xi_{13}}{18}; \\
	\Con[1, 3, 4, 4] &= \Con[3, 1, 4, 4] = \Con[4, 4, 1, 3] = \Con[4, 4, 3, 1] = -\frac{\zeta_{13}}{18} + \frac{\xi_{13}}{18}; \\
	\Con[1, 4, 2, 2] &= \Con[2, 2, 1, 4] = \Con[2, 2, 4, 1] = \Con[4, 1, 2, 2] = -\frac{\zeta_{12}}{18} + \frac{\xi_{12}}{18}; \\
	\Con[1, 4, 3, 3] &= \Con[3, 3, 1, 4] = \Con[3, 3, 4, 1] = \Con[4, 1, 3, 3] = \frac{\zeta_{12}}{18} - \frac{\xi_{12}}{18}; \\
	\Con[2, 2, 3, 4] &= \Con[2, 2, 4, 3] = \Con[3, 4, 2, 2] = \Con[4, 3, 2, 2] = \frac{\zeta_{23}}{18} + \frac{\xi_{23}}{18}; \\
	\Con[2, 3, 4, 4] &= \Con[3, 2, 4, 4] = \Con[4, 4, 2, 3] = \Con[4, 4, 3, 2] = \frac{\zeta_{12}}{18} + \frac{\xi_{12}}{18}; \\
	\Con[2, 4, 3, 3] &= \Con[3, 3, 2, 4] = \Con[3, 3, 4, 2] = \Con[4, 2, 3, 3] = \frac{\zeta_{13}}{18} + \frac{\xi_{13}}{18};
\end{align}
\subsection{Type VI}
\begin{align}
	\Con[1, 2, 1, 3] &= \Con[1, 3, 1, 2] = \Con[2, 1, 3, 1] = \Con[3, 1, 2, 1]; \\
	\Con[1, 2, 1, 4] &= \Con[1, 4, 1, 2] = \Con[2, 1, 4, 1] = \Con[4, 1, 2, 1]; \\
	\Con[1, 2, 2, 3] &= \Con[2, 1, 3, 2] = \Con[2, 3, 1, 2] = \Con[3, 2, 2, 1]; \\
	\Con[1, 2, 2, 4] &= \Con[2, 1, 4, 2] = \Con[2, 4, 1, 2] = \Con[4, 2, 2, 1]; \\
	\Con[1, 3, 1, 4] &= \Con[1, 4, 1, 3] = \Con[3, 1, 4, 1] = \Con[4, 1, 3, 1]; \\
	\Con[1, 3, 2, 3] &= \Con[2, 3, 1, 3] = \Con[3, 1, 3, 2] = \Con[3, 2, 3, 1]; \\
	\Con[1, 3, 3, 4] &= \Con[3, 1, 4, 3] = \Con[3, 4, 1, 3] = \Con[4, 3, 3, 1]; \\
	\Con[1, 4, 2, 4] &= \Con[2, 4, 1, 4] = \Con[4, 1, 4, 2] = \Con[4, 2, 4, 1]; \\
	\Con[1, 4, 3, 4] &= \Con[3, 4, 1, 4] = \Con[4, 1, 4, 3] = \Con[4, 3, 4, 1]; \\
	\Con[2, 3, 2, 4] &= \Con[2, 4, 2, 3] = \Con[3, 2, 4, 2] = \Con[4, 2, 3, 2]; \\
	\Con[2, 3, 3, 4] &= \Con[3, 2, 4, 3] = \Con[3, 4, 2, 3] = \Con[4, 3, 3, 2]; \\
	\Con[2, 4, 3, 4] &= \Con[3, 4, 2, 4] = \Con[4, 2, 4, 3] = \Con[4, 3, 4, 2]; \\
	\Con[1, 2, 3, 1] &= \Con[1, 3, 2, 1] = \Con[2, 1, 1, 3] = \Con[3, 1, 1, 2] \nonumber \\
	&= \frac{\zeta_{12}}{18} + \frac{\xi_{12}}{18} - \Con[1, 2, 1, 3]; \\
	\Con[1, 2, 4, 1] &= \Con[1, 4, 2, 1] = \Con[2, 1, 1, 4] = \Con[4, 1, 1, 2] \nonumber \\
	&= \frac{\zeta_{13}}{18} + \frac{\xi_{13}}{18} - \Con[1, 2, 1, 4]; \\
	\Con[1, 2, 4, 2] &= \Con[2, 1, 2, 4] = \Con[2, 4, 2, 1] = \Con[4, 2, 1, 2] \nonumber \\
	&= \frac{\zeta_{12}}{18} - \frac{\xi_{12}}{18} - \Con[1, 2, 2, 4]; \\
	\Con[1, 3, 3, 2] &= \Con[2, 3, 3, 1] = \Con[3, 1, 2, 3] = \Con[3, 2, 1, 3] \nonumber \\
	&= \frac{\zeta_{23}}{18} - \frac{\xi_{23}}{18} - \Con[1, 3, 2, 3]; \\
	\Con[1, 3, 4, 1] &= \Con[1, 4, 3, 1] = \Con[3, 1, 1, 4] = \Con[4, 1, 1, 3] \nonumber \\
	&= \frac{\zeta_{23}}{18} + \frac{\xi_{23}}{18} - \Con[1, 3, 1, 4]; \\
	\Con[1, 4, 4, 3] &= \Con[3, 4, 4, 1] = \Con[4, 1, 3, 4] = \Con[4, 3, 1, 4] \nonumber \\
	&= \frac{\zeta_{13}}{18} - \frac{\xi_{13}}{18} - \Con[1, 4, 3, 4]; \\
	\Con[1, 2, 3, 2] &= \Con[2, 1, 2, 3] = \Con[2, 3, 2, 1] = \Con[3, 2, 1, 2] \nonumber \\
	&= -\frac{\zeta_{13}}{18} + \frac{\xi_{13}}{18} - \Con[1, 2, 2, 3]; \\
	\Con[1, 3, 4, 3] &= \Con[3, 1, 3, 4] = \Con[3, 4, 3, 1] = \Con[4, 3, 1, 3] \nonumber \\
	&= -\frac{\zeta_{12}}{18} + \frac{\xi_{12}}{18} - \Con[1, 3, 3, 4]; \\
	\Con[1, 4, 4, 2] &= \Con[2, 4, 4, 1] = \Con[4, 1, 2, 4] = \Con[4, 2, 1, 4] \nonumber \\
	&= -\frac{\zeta_{23}}{18} + \frac{\xi_{23}}{18} - \Con[1, 4, 2, 4]; \\
	\Con[2, 3, 4, 2] &= \Con[2, 4, 3, 2] = \Con[3, 2, 2, 4] = \Con[4, 2, 2, 3] \nonumber \\
	&= -\frac{\zeta_{23}}{18} - \frac{\xi_{23}}{18} - \Con[2, 3, 2, 4]; \\
	\Con[2, 3, 4, 3] &= \Con[3, 2, 3, 4] = \Con[3, 4, 3, 2] = \Con[4, 3, 2, 3] \nonumber \\
	&= -\frac{\zeta_{13}}{18} - \frac{\xi_{13}}{18} - \Con[2, 3, 3, 4]; \\
	\Con[2, 4, 4, 3] &= \Con[3, 4, 4, 2] = \Con[4, 2, 3, 4] = \Con[4, 3, 2, 4] \nonumber \\
	&= -\frac{\zeta_{12}}{18} - \frac{\xi_{12}}{18} - \Con[2, 4, 3, 4];
\end{align}
\subsection{Type VII}
\begin{align}
	\Con[1, 2, 3, 4] &= \Con[2, 1, 4, 3] = \Con[3, 4, 1, 2] = \Con[4, 3, 2, 1]; \\
	\Con[1, 3, 2, 4] &= \Con[2, 4, 1, 3] = \Con[3, 1, 4, 2] = \Con[4, 2, 3, 1]; \\
	\Con[1, 4, 2, 3] &= \Con[2, 3, 1, 4] = \Con[3, 2, 4, 1] = \Con[4, 1, 3, 2]; \\
	\Con[1, 2, 4, 3] &= \Con[2, 1, 3, 4] = \Con[3, 4, 2, 1] = \Con[4, 3, 1, 2] \nonumber \\
	&= \frac{\zeta_{22}}{6} - \frac{\zeta_{33}}{6} - \frac{\xi_{22}}{6} + \frac{\xi_{33}}{6} - \Con[1, 2, 3, 4]; \\
	\Con[1, 4, 3, 2] &= \Con[2, 3, 4, 1] = \Con[3, 2, 1, 4] = \Con[4, 1, 2, 3] \nonumber \\
	&= \frac{\zeta_{11}}{6} - \frac{\zeta_{22}}{6} - \frac{\xi_{11}}{6} + \frac{\xi_{22}}{6} - \Con[1, 4, 2, 3]; \\
	\Con[1, 3, 4, 2] &= \Con[2, 4, 3, 1] = \Con[3, 1, 2, 4] = \Con[4, 2, 1, 3] \nonumber \\
	&= -\frac{\zeta_{11}}{6} + \frac{\zeta_{33}}{6} + \frac{\xi_{11}}{6} - \frac{\xi_{33}}{6} - \Con[1, 3, 2, 4].
\end{align}
	
	\appendix
	\section{The integration operator via Mathematica}
	In section 3, our integrands
	\begin{align}
	\frac{\p}{\p x^j}(\tilde{\eta}_{ii})+\sum_{k=1}^4h_{kj}\frac{\p}{\p x^j}((x^i)^2).
	\end{align}
	are of the type $r^{-6}P_3(x^1,x^2,x^3,x^4)$, where $P_3(x^1,x^2,x^3,x^4)$ denotes a homogeneous polynomial of $x_1,\cdots,x_4$ with degree 3, and $r^2=\sum_{i=1}^4 (x^i)^2.$ The terms that contribute non-zero integrals are similar to the following types (WLOG, consider the integral on $\Sigma^{\pm}_1$):
	\begin{align}
		&\frac{(x^1)^3}{r^6}dx^2\wedge dx^3\wedge dx^4,\quad\frac{x^1(x^2)^2}{r^6}dx^2\wedge dx^3\wedge dx^4,\\
		&\text{integrated on}\ \Sigma^{\pm}_1.
	\end{align}
	And in section 5, our integrands in these expressions are of the type $r^{-8}P_5(x^1,x^2,x^3,x^4)$. The terms with nonzero integrals are of the following types (WLOG, consider the integral on $\Sigma^{\pm}_1$):
	\begin{align}
		&\frac{x^1(x^2)^4}{r^8}dx^2\wedge dx^3\wedge dx^4,\quad\frac{x^1(x^2)^2(x^3)^2}{r^8}dx^2\wedge dx^3\wedge dx^4,\\
		&\frac{(x^1)^3(x^2)^2}{r^8}dx^2\wedge dx^3\wedge dx^4,\quad\frac{(x^1)^5}{r^8}dx^2\wedge dx^3\wedge dx^4,\\
		&\text{integrated on}\ \Sigma^{\pm}_1.
	\end{align}
	These terms are integrated over$\Sigma^{+}_1$. Moreover, the exterior normal vector fields of $\Sigma^{-}_1$ have opposite orientations. Thus, after setting $x^1=\pm1$ and integrating respectively, one will get the same result. For example,
	\begin{align}
		\int_{\Sigma^{+}_1}\frac{x^1(x^2)^4}{r^8}dx^2\wedge dx^3\wedge dx^4=\int_{\Sigma^{-}_1}\frac{x^1(x^2)^4}{r^8}dx^2\wedge dx^3\wedge dx^4=\int_{[-1,1]^3}\frac{(x^2)^4}{r^8}dx^2dx^3dx^4.
	\end{align}
	Therefore, we denote our integrals as (here $r^2=1+(x^2)^2+(x^3)^2+(x^4)^2$)
	\begin{align}
		c_1&=\int_{[-1,1]^3}\frac{1}{r^6}dx^2dx^3dx^4,\\
		c_2&=\int_{[-1,1]^3}\frac{(x^2)^2}{r^6}dx^2dx^3dx^4,\\
		\sigma_1&=\int_{[-1,1]^3}\frac{(x^2)^4}{r^8}dx^2dx^3dx^4,\\
		\sigma_2&=\int_{[-1,1]^3}\frac{(x^2)^2(x^3)^2}{r^8}dx^2dx^3dx^4,\\
		\sigma_3&=\int_{[-1,1]^3}\frac{(x^2)^2}{r^8}dx^2dx^3dx^4,\\
		\sigma_4&=\int_{[-1,1]^3}\frac{1}{r^8}dx^2dx^3dx^4,
	\end{align}
	We note that
	\begin{align}
		8 (3 \sigma_1+6 \sigma_2+6 \sigma_3+\sigma_4)=2\pi^2.
	\end{align}
	To show this, we have
	\begin{align}
		&(3 \sigma_1+6 \sigma_2+6 \sigma_3+\sigma_4)\\
		&=\int_{[-1,1]^3}\left.\frac{3x^1(x^2)^4+6x^1(x^2)^2(x^3)^2+6(x^1)^3(x^2)^2+(x^1)^5}{r^8}\right|_{x^1=1}dx^2dx^3dx^4\\
		&=\int_{[-1,1]^3}\frac{1}{r^8}\left(x^1(\sum_{i\neq2}(x^i)^4+\sum_{i< j,i\neq2,j\neq2}2(x^i)^2(x^j)^2\right.\\
		&\left.\left.+\sum_{i\neq2}2(x^i)^2(x^1)^4+(x^1)^4)\right)\right|_{x^1=1}dx^2dx^3dx^4\\
		&=\int_{[-1,1]^3}\left.\frac{x^1}{r^4}\right|_{x^1=1}dx^2dx^3dx^4.
	\end{align}
	Thus
	\begin{align}
		&8 (3 \sigma_1+6 \sigma_2+6 \sigma_3+\sigma_4)\\
		&=8\int_{[-1,1]^3}\left.\frac{x^1}{r^4}\right|_{x^1=1}dx^2dx^3dx^4\\
		&=\sum_{j=1}^4\int_{\Sigma_j^+}\frac{x^j}{r^4}dx^1\cdots\widehat{dx^j}\cdots dx^4-\sum_{j=1}^4\int_{\Sigma_j^-}\frac{x^j}{r^4}dx^1\cdots\widehat{dx^j}\cdots dx^4\\
		&=-\frac{1}{2}\int_{\p[-1,1]^4}*_{\R^4}d(\frac{1}{r^2})\\
		&=\frac{1}{2}\int_{[-1,1]^4}\Delta_{\R^4}(\frac{1}{r^2})\\
		&=\frac{1}{2}\int_{B_1(0)}\Delta_{\R^4}(\frac{1}{r^2})\\
		&=\frac{1}{2}\int_{\mathbb{S}^3}-\frac{\p}{\p r}(\frac{1}{r^2})dS\\
		&=\int_{\mathbb{S}^3}dS=|\mathbb{S}^3|=2\pi^2.
	\end{align}
	In addition, we also have
	\[
	3 \sigma_1-6 \sigma_2-6 \sigma_3+\sigma_4=0.
	\]
	To see this, it suffices to show that $8(3\sigma_1+\sigma_4)=\pi^2$. In fact,
	\begin{align}
		&8(3\sigma_1+\sigma_4)\\
		&=2\sum_{k=1}^{4}\int_{[-1,1]^3}\frac{1+3(x^2)^4}{(1+(x^2)^2+(x^3)^2+(x^4)^2)^4}dx^2dx^3dx^4\\
		&=\int_{\p[-1,1]^4}\sum_{k=1}^{4}\frac{x^k((x^1)^4+(x^2)^4+(x^3)^4+(x^4)^4)}{r^8}\gamma_k,
	\end{align}
	while
	\begin{align}
		&d\left(\sum_{k=1}^{4}\frac{x^k((x^1)^4+(x^2)^4+(x^3)^4+(x^4)^4)}{r^8}\gamma_k\right)\\
		&=\left(\sum_{k=1}^{4}\frac{\p}{\p x^k}(\frac{x^k((x^1)^4+(x^2)^4+(x^3)^4+(x^4)^4)}{r^8})\right)dx^1\wedge dx^2\wedge dx^3\wedge dx^4\\
		&=0.
	\end{align}
	Thus
	\begin{align}
		&\int_{\p[-1,1]^4}\sum_{k=1}^{4}\frac{x^k((x^1)^4+(x^2)^4+(x^3)^4+(x^4)^4)}{r^8}\gamma_k\\
		&=\int_{\mathbb{S}^3}\sum_{k=1}^{4}\frac{x^k((x^1)^4+(x^2)^4+(x^3)^4+(x^4)^4)}{r^8}\gamma_k\\
		&=\int_{0}^{\pi}d\theta_1\int_{0}^{\pi}d\theta_2\int_{0}^{2\pi}Q(\theta_1,\theta_2,\theta_3)d\theta_3
		=\pi^2,
	\end{align}
	where we have used the sphere coordinate
	\[
	\left\{
	\begin{array}{rcl}
		x_1 & = & \cos \theta_1, \\
		x_2 & = & \sin \theta_1 \cos \theta_2, \\
		x_3 & = & \sin \theta_1 \sin \theta_2 \cos \theta_3, \\
		x_4 & = & \sin \theta_1 \sin \theta_2 \sin \theta_3,
	\end{array}
	\right.
	\]
	in which
	\[
	\theta_1 \in (0, \pi), \quad
	\theta_2 \in (0, \pi), \quad
	\theta_3 \in (0, 2\pi),
	\]
	and
	\begin{align}
		&Q(\theta_1,\theta_2,\theta_3)\\
		&=\sin ^2(\theta_1) \sin (\theta_2) \left(\sin ^4(\theta_1) \left(\frac{1}{4} \sin ^4(\theta_2) (\cos (4 \theta_3)+3)+\cos ^4(\theta_2)\right)+\cos ^4(\theta_1)\right).
	\end{align}
	The complete computational code can be found in Appendix H.\\
	Here we employ the following method to calculate the integrals with Mathematica efficiently:
	\begin{enumerate}
		\item Input: $r^{-6}P_3(x_1,x_2,x_3,x_4)$ (resp. $r^{-8}P_5(x_1,x_2,x_3,x_4)$).
		\item Multiply the input by $r^{6}$ (resp. $r^{8}$).
		\item If the integration takes place on $\Sigma^{\pm}_k$, set $x_k$ to be $\pm 1$. $(k=1,2,3,4)$
		\item Construct an operator that can map the polynomial on the numerator to $c_i$ or 0 (resp. $\sigma_i$ or 0).
		\item Apply the operator to the polynomial.
		\item Output: Linear combination of $c_i$ (resp. $\sigma_i$).
	\end{enumerate}
	Let $p_n(x_1,x_2,x_3,x_4)$ denote a polynomial of order $n$ $(n=3,5)$.\\
	For $n=3$, the operator is defined as
	\begin{align}
		F(p_3):=\sum_{i,j,k,l=0}^1 \frac{1}{16}p_3\left((-1)^i,(-1)^j,(-1)^k,(-1)^l\right)c_1+p_3(0,0,0,0)(c_2-c_1)
	\end{align}
	This operator will map 1 to $c_1$, $(x^i)^2$ into $c_2$, and others to 0.\\
	For $n=5$, the operator is defined as
	\begin{align}
		G(p_5):=&G_1(p_5)\left(-\frac{\sigma _1}{12}+\frac{\sigma _2}{16}+\frac{\sigma _3}{12}\right)+G_2(p_5) \left(-\frac{11 \sigma _1}{4}+3 \sigma _2-\frac{5 \sigma _3}{4}+\sigma _4\right)\\
		&+G_3(p_5) \left(\frac{\sigma _1}{192}-\frac{\sigma _3}{192}\right)+G_4(p_5) \left(\frac{\sigma _1}{2}-\frac{\sigma _2}{2}\right),
	\end{align}
	where
	\begin{align}
		G_1(p_5)&:=\sum_{i,j,k,l=0}^1 p_5\left((-1)^i,(-1)^j,(-1)^k,(-1)^l\right),\\
		G_2(p_5)&:=p_5(0,0,0,0),\\
		G_3(p_5)&:=\sum_{i,j,k,l=0}^1 p_5\left(2\cdot(-1)^i,2\cdot(-1)^j,2\cdot(-1)^k,2\cdot(-1)^l\right),\\
		G_4(p_5)&:=\sum_{i=0}^{1}p_5\left((-1)^i,0,0,0\right)+p_5\left(0,(-1)^i,0,0\right)\\
		&+p_5\left(0,0,(-1)^i,0\right)+p_5\left(0,0,0,(-1)^i\right).\\
	\end{align}
	This operator will map $(x^i)^4$ to $\sigma_1$, $(x^i)^2(x^j)^2\ (i\neq j)$ to $\sigma_2$, $(x^i)^2$ to $\sigma_3$, 1 to $\sigma_4$, and others to 0.
	\subsection{Integration operator, n=3}
	The Mathematica code is as follows.
	\begin{lstlisting}
f[inp_] := 
Sum[(inp /. {x[1] -> (-1)^i, x[2] -> (-1)^j, x[3] -> (-1)^k, 
	x[4] -> (-1)^l}), {i, 0, 1}, {j, 0, 1}, {k, 0, 1}, {l, 0, 
	1}] \[Sigma][1]/
16 + (inp /. {x[1] -> 0, x[2] -> 0, x[3] -> 0, 
	x[4] -> 0}) (\[Sigma][2] - \[Sigma][1]) // Expand;
	\end{lstlisting}
	
	\subsection{Integration operator, n=5}
	The Mathematica code is as follows.
	\begin{lstlisting}
Process4E5TPolynomial := 
Function[{iinp}, (a[
1]*(-(\[Sigma][1]/12) + \[Sigma][2]/16 + \[Sigma][3]/12) + 
a[2]*(-((11 \[Sigma][1])/4) + 
3 \[Sigma][2] - (5 \[Sigma][3])/4 + \[Sigma][4]) + 
a[3]*(\[Sigma][1]/192 - \[Sigma][3]/192) + 
a[4]*(\[Sigma][1]/2 - \[Sigma][2]/2)) /. {a[1] -> 
	Sum[(iinp /. {x[1] -> (-1)^i, x[2] -> (-1)^j, x[3] -> (-1)^k, 
		x[4] -> (-1)^l}), {i, 0, 1}, {j, 0, 1}, {k, 0, 1}, {l, 0, 
		1}], a[2] -> (iinp /. {x[1] -> 0, x[2] -> 0, x[3] -> 0, 
		x[4] -> 0}), 
	a[3] -> Sum[(iinp /. {x[1] -> 2*(-1)^i, x[2] -> 2*(-1)^j, 
		x[3] -> 2*(-1)^k, x[4] -> 2*(-1)^l}), {i, 0, 1}, {j, 0, 
		1}, {k, 0, 1}, {l, 0, 1}], 
	a[4] -> Sum[
	iinp /. {x[1] -> gener[i], x[2] -> gener[i + 1], 
		x[3] -> gener[i + 2], x[4] -> gener[i + 3]}, {i, 0, 7}] /. 
	gener[igenerProcess4E5TPolynomial_] -> (Cos[
	Pi*igenerProcess4E5TPolynomial] + 
	Abs[Cos[Pi*igenerProcess4E5TPolynomial]])*
	Cos[Pi/4*igenerProcess4E5TPolynomial]/2} // Expand];
	\end{lstlisting}
	\section{Harmonic Function.nb}
	\begin{lstlisting}
x[1] = x1;
x[2] = x2;
x[3] = x3;
x[4] = x4;
dx[1] = dx1;
dx[2] = dx2;
dx[3] = dx3;
dx[4] = dx4;
X[1] = \[Zeta]11;
X[2] = \[Zeta]12;
X[3] = \[Zeta]13;
X[4] = \[Zeta]22;
X[5] = \[Zeta]23;
X[6] = \[Zeta]33;
X[7] = \[Xi]11;
X[8] = \[Xi]12;
X[9] = \[Xi]13;
X[10] = \[Xi]22;
X[11] = \[Xi]23;
X[12] = \[Xi]33;
r = Sqrt[Sum[x[i]^2, {i, 1, 4}]];
J[1] = {{0, -1, 0, 0}, {1, 0, 0, 0}, {0, 0, 0, -1}, {0, 0, 1, 0}};
J[2] = {{0, 0, -1, 0}, {0, 0, 0, 1}, {1, 0, 0, 0}, {0, -1, 0, 0}};
J[3] = J[1] . J[2];
\[Alpha][0] = {{x1}, {x2}, {x3}, {x4}};
Do[\[Alpha][i] = Transpose[J[i]] . \[Alpha][0], {i, 1, 4}];
Do[A[i, i] = \[Alpha][i] . Transpose[\[Alpha][i]], {i, 1, 3}];
Do[If[i != j, A[i, j] = \[Alpha][i] . Transpose[\[Alpha][j]]], {i, 1, 
	3}, {j, 1, 3}];
hp = (-(3/2) r^6)^(-1)*(\[Zeta]11 (2 A[1, 1] - A[2, 2] - 
A[3, 3]) + \[Zeta]22 (2 A[2, 2] - A[1, 1] - 
A[3, 3]) + \[Zeta]33 (2 A[3, 3] - A[1, 1] - 
A[2, 2]) + \[Zeta]12 (A[1, 2] + 
A[2, 1]) + \[Zeta]13 (A[1, 3] + 
A[3, 1]) + \[Zeta]23 (A[2, 3] + A[3, 2]));

R = {{1, 0, 0, 0}, {0, -1, 0, 0}, {0, 0, -1, 0}, {0, 0, 0, -1}};
Do[m\[Alpha][i] = 
Transpose[R] . Transpose[J[i]] . R . \[Alpha][0], {i, 1, 4}];
Do[mA[i, i] = m\[Alpha][i] . Transpose[m\[Alpha][i]], {i, 1, 3}];
Do[If[i != j, mA[i, j] = m\[Alpha][i] . Transpose[m\[Alpha][j]]], {i, 
	1, 3}, {j, 1, 3}];
hm = (-(3/2) r^6)^(-1)*(\[Xi]11 (2 mA[1, 1] - mA[2, 2] - 
mA[3, 3]) + \[Xi]22 (2 mA[2, 2] - mA[1, 1] - 
mA[3, 3]) + \[Xi]33 (2 mA[3, 3] - mA[1, 1] - 
mA[2, 2]) + \[Xi]12 (mA[1, 2] + 
mA[2, 1]) + \[Xi]13 (mA[1, 3] + 
mA[3, 1]) + \[Xi]23 (mA[2, 3] + mA[3, 2]));
h = hp + hm;
Do[t\[Eta][i, j] = h[[i, j]]*r^2/4, {i, 1, 4}, {j, 1, 4}];
Do[Q[i, j, l] = 
D[t\[Eta][i, j], x[l]] + 
Sum[h[[k, l]]*D[x[i]*x[j], x[k]], {k, 1, 4}], {i, 1, 4}, {j, 1, 
	4}, {l, 1, 4}];
Do[QQ[i, j, l] = r^6*Q[i, j, l], {i, 1, 4}, {j, 1, 4}, {l, 1, 4}];
Capp[inp_] := 
Sum[(inp /. {x[1] -> (-1)^i, x[2] -> (-1)^j, x[3] -> (-1)^k, 
	x[4] -> (-1)^l}), {i, 0, 1}, {j, 0, 1}, {k, 0, 1}, {l, 0, 
	1}] \[Sigma][1]/
16 + (inp /. {x[1] -> 0, x[2] -> 0, x[3] -> 0, 
	x[4] -> 0}) (\[Sigma][2] - \[Sigma][1]) // Expand;
Do[Int[i, j] = 
Sum[Capp[(QQ[i, j, l] /. x[l] -> 1) - (QQ[i, j, l] /. 
x[l] -> -1)], {l, 1, 4}], {i, 1, 4}, {j, 1, 4}];
Do[Print[Int[i, j]], {i, 1, 4}, {j, 1, 4}];
Do[Int[i, j] = Simplify[Int[i, j]], {i, 1, 4}, {j, 1, 4}];
Do[Print[Int[i, j]], {i, 1, 4}, {j, 1, 4}];		
	\end{lstlisting}
	\section{Differential of Harmonic Functions.nb}
	\begin{lstlisting}
x[1] = x1;
x[2] = x2;
x[3] = x3;
x[4] = x4;
dx[1] = dx1;
dx[2] = dx2;
dx[3] = dx3;
dx[4] = dx4;
X[1] = \[Zeta]11;
X[2] = \[Zeta]12;
X[3] = \[Zeta]13;
X[4] = \[Zeta]22;
X[5] = \[Zeta]23;
X[6] = \[Zeta]33;
X[7] = \[Xi]11;
X[8] = \[Xi]12;
X[9] = \[Xi]13;
X[10] = \[Xi]22;
X[11] = \[Xi]23;
X[12] = \[Xi]33;
r = Sqrt[Sum[x[i]^2, {i, 1, 4}]];
J[1] = {{0, -1, 0, 0}, {1, 0, 0, 0}, {0, 0, 0, -1}, {0, 0, 1, 0}};
J[2] = {{0, 0, -1, 0}, {0, 0, 0, 1}, {1, 0, 0, 0}, {0, -1, 0, 0}};
J[3] = J[1] . J[2];
\[Alpha][0] = {{x1}, {x2}, {x3}, {x4}};
Do[\[Alpha][i] = Transpose[J[i]] . \[Alpha][0], {i, 1, 4}];
Do[A[i, i] = \[Alpha][i] . Transpose[\[Alpha][i]], {i, 1, 3}];
Do[If[i != j, A[i, j] = \[Alpha][i] . Transpose[\[Alpha][j]]], {i, 1, 
	3}, {j, 1, 3}];
hp = (-(3/2) r^6)^(-1)*(\[Zeta]11 (2 A[1, 1] - A[2, 2] - 
A[3, 3]) + \[Zeta]22 (2 A[2, 2] - A[1, 1] - 
A[3, 3]) + \[Zeta]33 (2 A[3, 3] - A[1, 1] - 
A[2, 2]) + \[Zeta]12 (A[1, 2] + 
A[2, 1]) + \[Zeta]13 (A[1, 3] + 
A[3, 1]) + \[Zeta]23 (A[2, 3] + A[3, 2]));

R = {{1, 0, 0, 0}, {0, -1, 0, 0}, {0, 0, -1, 0}, {0, 0, 0, -1}};
Do[m\[Alpha][i] = 
Transpose[R] . Transpose[J[i]] . R . \[Alpha][0], {i, 1, 4}];
Do[mA[i, i] = m\[Alpha][i] . Transpose[m\[Alpha][i]], {i, 1, 3}];
Do[If[i != j, mA[i, j] = m\[Alpha][i] . Transpose[m\[Alpha][j]]], {i, 
	1, 3}, {j, 1, 3}];
hm = (-(3/2) r^6)^(-1)*(\[Xi]11 (2 mA[1, 1] - mA[2, 2] - 
mA[3, 3]) + \[Xi]22 (2 mA[2, 2] - mA[1, 1] - 
mA[3, 3]) + \[Xi]33 (2 mA[3, 3] - mA[1, 1] - 
mA[2, 2]) + \[Xi]12 (mA[1, 2] + 
mA[2, 1]) + \[Xi]13 (mA[1, 3] + 
mA[3, 1]) + \[Xi]23 (mA[2, 3] + mA[3, 2]));
h = hp + hm;
Do[t\[Eta][i, j] = h[[i, j]]*r^2/4, {i, 1, 4}, {j, 1, 4}];
Do[B[1, i1, i2] = 
Sum[D[h[[i1, i2]] + Sum[x[i1]*D[h[[i2, j]], x[j]], {j, 1, 4}], 
x[k]]*dx[k], {k, 1, 4}], {i1, 1, 4}, {i2, 1, 4}];
Do[B[2, i1, i2] = 
Sum[(Boole[i1 == l && i2 == k] - Boole[i1 == k && i2 == l])*
D[Boole[i == k]*h[[j, l]] + h[[i, k]]*Boole[j == l], x[j]]*
dx[i], {i, 1, 4}, {j, 1, 4}, {k, 1, 4}, {l, 1, 4}], {i1, 1, 
	4}, {i2, 1, 4}];

Do[L[i1, i2] = B[1, i1, i2] + B[2, i1, i2], {i1, 1, 4}, {i2, 1, 4}];
Do[t\[Mu][i1, i2] = L[i1, i2]*r^2/(12), {i1, 1, 4}, {i2, 1, 4}];
Do[t\[Omega][i1, i2] = 
t\[Mu][i1, i2] + 
Sum[(Con[i1, i2, k, l]*x[k])/r^4*dx[l], {k, 1, 4}, {l, 1, 
	4}], {i1, 1, 4}, {i2, 1, 4}];
Do[\[Omega][i1, i2] = x[i1]*dx[i2] - t\[Omega][i1, i2], {i1, 1, 
	4}, {i2, 1, 4}];
Do[u[i, i] = x[i]*x[i] - t\[Eta][i, i] + CVol/r^2, {i, 1, 4}];
Do[If[i != j, u[i, j] = x[i]*x[j] - t\[Eta][i, j]], {i, 1, 4}, {j, 1, 
	4}];
Do[du[i, j] = Sum[D[u[i, j], x[k]]*dx[k], {k, 1, 4}], {i, 1, 4}, {j, 
	1, 4}];
Do[Dif[i, j] = 
Simplify[(du[i, j] - \[Omega][i, j] - \[Omega][j, i])]*9 r^4, {i, 
	1, 4}, {j, 1, 4}];
Do[DDD[i, j, k, l] = D[Dif[i, j], x[k], dx[l]], {i, 1, 4}, {j, 1, 
	4}, {k, 1, 4}, {l, 1, 4}];
Do[Print[DDD[i, j, k, l]], {i, 1, 4}, {j, 1, 4}, {k, 1, 4}, {l, 1, 
	4}];		
	\end{lstlisting}	
	\section{Divergence Arguments.nb}
	\begin{lstlisting}
x[1] = x1;
x[2] = x2;
x[3] = x3;
x[4] = x4;
dx[1] = dx1;
dx[2] = dx2;
dx[3] = dx3;
dx[4] = dx4;
X[1] = \[Zeta]11;
X[2] = \[Zeta]12;
X[3] = \[Zeta]13;
X[4] = \[Zeta]22;
X[5] = \[Zeta]23;
X[6] = \[Zeta]33;
X[7] = \[Xi]11;
X[8] = \[Xi]12;
X[9] = \[Xi]13;
X[10] = \[Xi]22;
X[11] = \[Xi]23;
X[12] = \[Xi]33;
r = Sqrt[Sum[x[i]^2, {i, 1, 4}]];
J[1] = {{0, -1, 0, 0}, {1, 0, 0, 0}, {0, 0, 0, -1}, {0, 0, 1, 0}};
J[2] = {{0, 0, -1, 0}, {0, 0, 0, 1}, {1, 0, 0, 0}, {0, -1, 0, 0}};
J[3] = J[1] . J[2];
\[Alpha][0] = {{x1}, {x2}, {x3}, {x4}};
Do[\[Alpha][i] = Transpose[J[i]] . \[Alpha][0], {i, 1, 4}];
Do[A[i, i] = \[Alpha][i] . Transpose[\[Alpha][i]], {i, 1, 3}];
Do[If[i != j, A[i, j] = \[Alpha][i] . Transpose[\[Alpha][j]]], {i, 1, 
	3}, {j, 1, 3}];
hp = (-(3/2)
r^6)^(-1)*(\[Zeta]11 (2 A[1, 1] - A[2, 2] - 
A[3, 3]) + \[Zeta]22 (2 A[2, 2] - A[1, 1] - 
A[3, 3]) + \[Zeta]33 (2 A[3, 3] - A[1, 1] - 
A[2, 2]) + \[Zeta]12 (A[1, 2] + 
A[2, 1]) + \[Zeta]13 (A[1, 3] + 
A[3, 1]) + \[Zeta]23 (A[2, 3] + A[3, 2]));

R = {{1, 0, 0, 0}, {0, -1, 0, 0}, {0, 0, -1, 0}, {0, 0, 0, -1}};
Do[m\[Alpha][i] = 
Transpose[R] . Transpose[J[i]] . R . \[Alpha][0], {i, 1, 4}];
Do[mA[i, i] = m\[Alpha][i] . Transpose[m\[Alpha][i]], {i, 1, 3}];
Do[If[i != j, mA[i, j] = m\[Alpha][i] . Transpose[m\[Alpha][j]]], {i, 
	1, 3}, {j, 1, 3}];
hm = (-(3/2)
r^6)^(-1)*(\[Xi]11 (2 mA[1, 1] - mA[2, 2] - 
mA[3, 3]) + \[Xi]22 (2 mA[2, 2] - mA[1, 1] - 
mA[3, 3]) + \[Xi]33 (2 mA[3, 3] - mA[1, 1] - 
mA[2, 2]) + \[Xi]12 (mA[1, 2] + 
mA[2, 1]) + \[Xi]13 (mA[1, 3] + 
mA[3, 1]) + \[Xi]23 (mA[2, 3] + mA[3, 2]));
h = hp + hm;
Do[B[1, i1, i2] = 
Sum[D[h[[i1, i2]] + Sum[x[i1]*D[h[[i2, j]], x[j]], {j, 1, 4}], 
x[k]]*dx[k], {k, 1, 4}], {i1, 1, 4}, {i2, 1, 4}];

Do[B[2, i1, i2] = 
Sum[(Boole[i1 == l && i2 == k] - Boole[i1 == k && i2 == l])*
D[Boole[i == k]*h[[j, l]] + h[[i, k]]*Boole[j == l], x[j]]*
dx[i], {i, 1, 4}, {j, 1, 4}, {k, 1, 4}, {l, 1, 4}], {i1, 1, 
	4}, {i2, 1, 4}]; Do[
L[i1, i2] = B[1, i1, i2] + B[2, i1, i2], {i1, 1, 4}, {i2, 1, 4}];
Do[t\[Mu][i1, i2] = L[i1, i2]*r^2/(12), {i1, 1, 4}, {i2, 1, 4}]; Do[
t\[Omega][i1, i2] = 
t\[Mu][i1, i2] + 
Sum[(Con[i1, i2, k, l]*x[k])/r^4*dx[l], {k, 1, 4}, {l, 1, 4}], {i1,
	1, 4}, {i2, 1, 4}];
Do[\[Omega][i1, i2] = x[i1]*dx[i2] - t\[Omega][i1, i2], {i1, 1, 
	4}, {i2, 1, 4}];
Do[div[i1, i2] = 
h[[i1, i2]] + Sum[D[t\[Omega][i1, i2], dx[l], x[l]], {l, 1, 4}] + 
Sum[x[i1]*D[h[[i2, j]], x[j]], {j, 1, 4}], {i1, 1, 4}, {i2, 1, 4}];
Do[div[i1, i2] = Simplify[div[i1, i2]*9*r^6], {i1, 1, 4}, {i2, 1, 4}];
Array[Eq, {4, 4, 4, 4}];
Array[a, {4, 4, 4, 4}];
Do[Eq[i1, i2, j, k] = Simplify[D[div[i1, i2], x[j], x[k]]], {i1, 1, 
	4}, {i2, 1, 4}, {j, 1, 4}, {k, 1, 4}];
Do[Print[a[i1, i2, j, k] + Eq[i1, i2, j, k]], {i1, 1, 4}, {i2, 1, 
	4}, {j, 1, 4}, {k, 1, 4}];
Do[AAA[i1, i2, i3, i4] = Eq[i1, i2, i3, i4], {i1, 1, 4}, {i2, 1, 
	4}, {i3, 1, 4}, {i4, 1, 4}];
	\end{lstlisting}
	\section{Integral of the Laplacian.nb}
	\begin{lstlisting}
x[1] = x1;
x[2] = x2;
x[3] = x3;
x[4] = x4;
dx[1] = dx1;
dx[2] = dx2;
dx[3] = dx3;
dx[4] = dx4;
r = Sqrt[Sum[x[i]^2, {i, 1, 4}]];
J[1] = {{0, -1, 0, 0}, {1, 0, 0, 0}, {0, 0, 0, -1}, {0, 0, 1, 0}};
J[2] = {{0, 0, -1, 0}, {0, 0, 0, 1}, {1, 0, 0, 0}, {0, -1, 0, 0}};
J[3] = J[1] . J[2];
\[Alpha][0] = {{x1}, {x2}, {x3}, {x4}};
Do[\[Alpha][i] = Transpose[J[i]] . \[Alpha][0], {i, 1, 4}];
Do[A[i, i] = \[Alpha][i] . Transpose[\[Alpha][i]], {i, 1, 3}];
Do[If[i != j, A[i, j] = \[Alpha][i] . Transpose[\[Alpha][j]]], {i, 1, 
	3}, {j, 1, 3}];
hp = (-(3/2) r^6)^(-1)*(\[Zeta]11 (2 A[1, 1] - A[2, 2] - 
A[3, 3]) + \[Zeta]22 (2 A[2, 2] - A[1, 1] - 
A[3, 3]) + \[Zeta]33 (2 A[3, 3] - A[1, 1] - 
A[2, 2]) + \[Zeta]12 (A[1, 2] + 
A[2, 1]) + \[Zeta]13 (A[1, 3] + 
A[3, 1]) + \[Zeta]23 (A[2, 3] + A[3, 2]));

R = {{1, 0, 0, 0}, {0, -1, 0, 0}, {0, 0, -1, 0}, {0, 0, 0, -1}};
Do[m\[Alpha][i] = 
Transpose[R] . Transpose[J[i]] . R . \[Alpha][0], {i, 1, 4}];
Do[mA[i, i] = m\[Alpha][i] . Transpose[m\[Alpha][i]], {i, 1, 3}];
Do[If[i != j, mA[i, j] = m\[Alpha][i] . Transpose[m\[Alpha][j]]], {i, 
	1, 3}, {j, 1, 3}];
hm = (-(3/2) r^6)^(-1)*(\[Xi]11 (2 mA[1, 1] - mA[2, 2] - 
mA[3, 3]) + \[Xi]22 (2 mA[2, 2] - mA[1, 1] - 
mA[3, 3]) + \[Xi]33 (2 mA[3, 3] - mA[1, 1] - 
mA[2, 2]) + \[Xi]12 (mA[1, 2] + 
mA[2, 1]) + \[Xi]13 (mA[1, 3] + 
mA[3, 1]) + \[Xi]23 (mA[2, 3] + mA[3, 2]));
h = hp + hm;
Process4E5TPolynomial := 
Function[{iinp}, (a[
1]*(-(\[Sigma][1]/12) + \[Sigma][2]/16 + \[Sigma][3]/12) + 
a[2]*(-((11 \[Sigma][1])/4) + 
3 \[Sigma][2] - (5 \[Sigma][3])/4 + \[Sigma][4]) + 
a[3]*(\[Sigma][1]/192 - \[Sigma][3]/192) + 
a[4]*(\[Sigma][1]/2 - \[Sigma][2]/2)) /. {a[1] -> 
	Sum[(iinp /. {x[1] -> (-1)^i, x[2] -> (-1)^j, x[3] -> (-1)^k, 
		x[4] -> (-1)^l}), {i, 0, 1}, {j, 0, 1}, {k, 0, 1}, {l, 0, 
		1}], a[2] -> (iinp /. {x[1] -> 0, x[2] -> 0, x[3] -> 0, 
		x[4] -> 0}), 
	a[3] -> Sum[(iinp /. {x[1] -> 2*(-1)^i, x[2] -> 2*(-1)^j, 
		x[3] -> 2*(-1)^k, x[4] -> 2*(-1)^l}), {i, 0, 1}, {j, 0, 
		1}, {k, 0, 1}, {l, 0, 1}], 
	a[4] -> Sum[
	iinp /. {x[1] -> gener[i], x[2] -> gener[i + 1], 
		x[3] -> gener[i + 2], x[4] -> gener[i + 3]}, {i, 0, 7}] /. 
	gener[igenerProcess4E5TPolynomial_] -> (Cos[
	Pi*igenerProcess4E5TPolynomial] + 
	Abs[Cos[Pi*igenerProcess4E5TPolynomial]])*
	Cos[Pi/4*igenerProcess4E5TPolynomial]/2} // Expand];
Frigus := Process4E5TPolynomial;

Do[B[1, i1, i2] = 
Sum[D[h[[i1, i2]] + Sum[x[i1]*D[h[[i2, j]], x[j]], {j, 1, 4}], 
x[k]]*dx[k], {k, 1, 4}], {i1, 1, 4}, {i2, 1, 4}];
Do[B[2, i1, i2] = 
Sum[(Boole[i1 == l && i2 == k] - Boole[i1 == k && i2 == l])*
D[Boole[i == k]*h[[j, l]] + h[[i, k]]*Boole[j == l], x[j]]*
dx[i], {i, 1, 4}, {j, 1, 4}, {k, 1, 4}, {l, 1, 4}], {i1, 1, 
	4}, {i2, 1, 4}];

Do[L[i1, i2] = B[1, i1, i2] + B[2, i1, i2], {i1, 1, 4}, {i2, 1, 4}];
Do[t\[Mu][i1, i2] = L[i1, i2]*r^2/(12), {i1, 1, 4}, {i2, 1, 4}];
Do[t\[Omega][i1, i2] = 
t\[Mu][i1, i2] + 
Sum[(Con[i1, i2, k, l]*x[k])/r^4*dx[l], {k, 1, 4}, {l, 1, 
	4}], {i1, 1, 4}, {i2, 1, 4}];
Do[\[Omega][i1, i2] = x[i1]*dx[i2] - t\[Omega][i1, i2], {i1, 1, 
	4}, {i2, 1, 4}];
Do[\[CapitalOmega][i1, i2, i3, 
i4] = (-1) Sum[
Boole[(j3 < j4) && (j1 != j2 != j3 != j4)] Signature[{j1, j2, j3,
	j4}]*Signature[{j3, j4, t}]*dx[j3]*dx[j4]*
dx[t] (Boole[i1 == j1 && i2 == j2] D[t\[Omega][i3, i4], 
dx[t]] + (Boole[i1 == j1] h[[i2, j2]] + 
h[[i1, j1]] Boole[i2 == j2] + 
D[t\[Omega][i1, i2], dx[j2], x[j1]])*
x[i3] Boole[i4 == t]), {j1, 1, 4}, {j2, 1, 4}, {j3, 1, 
	4}, {j4, 1, 4}, {t, 1, 4}], {i1, 1, 4}, {i2, 1, 4}, {i3, 1, 
	4}, {i4, 1, 4}];
Do[\[CapitalOmega]A[i1, i2, i3, i4, p] = (-1)^(p - 1)*
D[\[CapitalOmega][i1, i2, i3, i4]*dx[p], dx1, dx2, dx3, dx4], {i1,
	1, 4}, {i2, 1, 4}, {i3, 1, 4}, {i4, 1, 4}, {p, 1, 4}];
Do[\[CapitalOmega]A[i1, i2, i3, i4, 
p] = \[CapitalOmega]A[i1, i2, i3, i4, p]*r^8, {i1, 1, 4}, {i2, 1, 
	4}, {i3, 1, 4}, {i4, 1, 4}, {p, 1, 4}];

Do[\[CapitalOmega]B[i1, i2, i3, i4, 
p] = (\[CapitalOmega]A[i1, i2, i3, i4, p] /. 
x[p] -> 1) - (\[CapitalOmega]A[i1, i2, i3, i4, p] /. 
x[p] -> -1), {i1, 1, 4}, {i2, 1, 4}, {i3, 1, 4}, {i4, 1, 4}, {p,
	1, 4}];
Do[Int2[i1, i2, i3, i4] = 
Sum[Frigus[\[CapitalOmega]B[i1, i2, i3, i4, p]], {p, 1, 4}], {i1, 
	1, 4}, {i2, 1, 4}, {i3, 1, 4}, {i4, 1, 4}];

Do[Print[a[i1, i2, i3, i4] + 
Simplify[(Int2[i1, i2, i3, i4] - 
Int2[i3, i4, i1, i2]), {3 \[Sigma][1] + 6 \[Sigma][2] + 
	6 \[Sigma][3] + \[Sigma][4] == Pi^2/4, 
	3 \[Sigma][1] - 6 \[Sigma][2] - 6 \[Sigma][3] + \[Sigma][4] == 
	0}]], {i1, 1, 4}, {i2, 1, 4}, {i3, 1, 4}, {i4, 1, 4}];
Do[BBB[i1, i2, i3, i4] = 
Simplify[(Int2[i1, i2, i3, i4] - 
Int2[i3, i4, i1, i2]), {3 \[Sigma][1] + 6 \[Sigma][2] + 
	6 \[Sigma][3] + \[Sigma][4] == Pi^2/4, 
	3 \[Sigma][1] - 6 \[Sigma][2] - 6 \[Sigma][3] + \[Sigma][4] == 
	0}], {i1, 1, 4}, {i2, 1, 4}, {i3, 1, 4}, {i4, 1, 4}];		
	\end{lstlisting}
	\section{Integral of the Covariant Derivative.nb}
	\begin{lstlisting}
x[1] = x1;
x[2] = x2;
x[3] = x3;
x[4] = x4;
dx[1] = dx1;
dx[2] = dx2;
dx[3] = dx3;
dx[4] = dx4;
r = Sqrt[Sum[x[i]^2, {i, 1, 4}]];
J[1] = {{0, -1, 0, 0}, {1, 0, 0, 0}, {0, 0, 0, -1}, {0, 0, 1, 0}};
J[2] = {{0, 0, -1, 0}, {0, 0, 0, 1}, {1, 0, 0, 0}, {0, -1, 0, 0}};
J[3] = J[1] . J[2];
\[Alpha][0] = {{x1}, {x2}, {x3}, {x4}};
Do[\[Alpha][i] = Transpose[J[i]] . \[Alpha][0], {i, 1, 4}];
Do[A[i, i] = \[Alpha][i] . Transpose[\[Alpha][i]], {i, 1, 3}];
Do[If[i != j, A[i, j] = \[Alpha][i] . Transpose[\[Alpha][j]]], {i, 1, 
	3}, {j, 1, 3}];
hp = (-(3/2) r^6)^(-1)*(\[Zeta]11 (2 A[1, 1] - A[2, 2] - 
A[3, 3]) + \[Zeta]22 (2 A[2, 2] - A[1, 1] - 
A[3, 3]) + \[Zeta]33 (2 A[3, 3] - A[1, 1] - 
A[2, 2]) + \[Zeta]12 (A[1, 2] + 
A[2, 1]) + \[Zeta]13 (A[1, 3] + 
A[3, 1]) + \[Zeta]23 (A[2, 3] + A[3, 2]));

R = {{1, 0, 0, 0}, {0, -1, 0, 0}, {0, 0, -1, 0}, {0, 0, 0, -1}};
Do[m\[Alpha][i] = 
Transpose[R] . Transpose[J[i]] . R . \[Alpha][0], {i, 1, 4}];
Do[mA[i, i] = m\[Alpha][i] . Transpose[m\[Alpha][i]], {i, 1, 3}];
Do[If[i != j, mA[i, j] = m\[Alpha][i] . Transpose[m\[Alpha][j]]], {i, 
	1, 3}, {j, 1, 3}];
hm = (-(3/2) r^6)^(-1)*(\[Xi]11 (2 mA[1, 1] - mA[2, 2] - 
mA[3, 3]) + \[Xi]22 (2 mA[2, 2] - mA[1, 1] - 
mA[3, 3]) + \[Xi]33 (2 mA[3, 3] - mA[1, 1] - 
mA[2, 2]) + \[Xi]12 (mA[1, 2] + 
mA[2, 1]) + \[Xi]13 (mA[1, 3] + 
mA[3, 1]) + \[Xi]23 (mA[2, 3] + mA[3, 2]));
h = hp + hm;
Process4E5TPolynomial := 
Function[{iinp}, (a[
1]*(-(\[Sigma][1]/12) + \[Sigma][2]/16 + \[Sigma][3]/12) + 
a[2]*(-((11 \[Sigma][1])/4) + 
3 \[Sigma][2] - (5 \[Sigma][3])/4 + \[Sigma][4]) + 
a[3]*(\[Sigma][1]/192 - \[Sigma][3]/192) + 
a[4]*(\[Sigma][1]/2 - \[Sigma][2]/2)) /. {a[1] -> 
	Sum[(iinp /. {x[1] -> (-1)^i, x[2] -> (-1)^j, x[3] -> (-1)^k, 
		x[4] -> (-1)^l}), {i, 0, 1}, {j, 0, 1}, {k, 0, 1}, {l, 0, 
		1}], a[2] -> (iinp /. {x[1] -> 0, x[2] -> 0, x[3] -> 0, 
		x[4] -> 0}), 
	a[3] -> Sum[(iinp /. {x[1] -> 2*(-1)^i, x[2] -> 2*(-1)^j, 
		x[3] -> 2*(-1)^k, x[4] -> 2*(-1)^l}), {i, 0, 1}, {j, 0, 
		1}, {k, 0, 1}, {l, 0, 1}], 
	a[4] -> Sum[
	iinp /. {x[1] -> gener[i], x[2] -> gener[i + 1], 
		x[3] -> gener[i + 2], x[4] -> gener[i + 3]}, {i, 0, 7}] /. 
	gener[igenerProcess4E5TPolynomial_] -> (Cos[
	Pi*igenerProcess4E5TPolynomial] + 
	Abs[Cos[Pi*igenerProcess4E5TPolynomial]])*
	Cos[Pi/4*igenerProcess4E5TPolynomial]/2} // Expand];
Frigus := Process4E5TPolynomial;
Array[B, {2, 4, 4}];
Do[B[1, i1, i2] = 
Sum[D[h[[i1, i2]] + Sum[x[i1]*D[h[[i2, j]], x[j]], {j, 1, 4}], 
x[k]]*dx[k], {k, 1, 4}], {i1, 1, 4}, {i2, 1, 4}];
Do[B[2, i1, i2] = 
Sum[(Boole[i1 == l && i2 == k] - Boole[i1 == k && i2 == l])*
D[Boole[i == k]*h[[j, l]] + h[[i, k]]*Boole[j == l], x[j]]*
dx[i], {i, 1, 4}, {j, 1, 4}, {k, 1, 4}, {l, 1, 4}], {i1, 1, 
	4}, {i2, 1, 4}];

Do[L[i1, i2] = B[1, i1, i2] + B[2, i1, i2], {i1, 1, 4}, {i2, 1, 4}];
Do[t\[Mu][i1, i2] = L[i1, i2]*r^2/(12), {i1, 1, 4}, {i2, 1, 4}];
Do[\[CapitalGamma][i, j, 
k] = (1/2) (D[h[[i, k]], x[j]] + D[h[[j, k]], x[i]] - 
D[h[[i, j]], x[k]]), {i, 1, 4}, {j, 1, 4}, {k, 1, 4}];
(*\[CapitalGamma][i,j,k] means \[CapitalGamma]_ij^k*)
Do[t\[Omega][i1, i2] = 
t\[Mu][i1, i2] + 
Sum[(Con[i1, i2, k, l]*x[k])/r^4*dx[l], {k, 1, 4}, {l, 1, 
	4}], {i1, 1, 4}, {i2, 1, 4}];
Do[\[Omega][i1, i2] = x[i1]*dx[i2] - t\[Omega][i1, i2], {i1, 1, 
	4}, {i2, 1, 4}];
Array[V, {4, 4, 4, 4, 4}];
Do[V[i1, i2, i3, i4, 
j] = (-1)*(D[t\[Omega][i1, i2], dx[i4]]*Boole[i3 == j] + 
x[i1]*(D[t\[Omega][i3, i4], x[j], dx[i2]] + 
x[i3]*\[CapitalGamma][i2, j, i4]) + 
x[i1]*Boole[i2 == i4]*h[[i3, j]] + 
x[i1]*Boole[i3 == j]*h[[i2, i4]]), {i1, 1, 4}, {i2, 1, 4}, {i3, 
	1, 4}, {i4, 1, 4}, {j, 1, 4}];
Do[VV[i1, i2, i3, i4, j] = r^8 V[i1, i2, i3, i4, j], {i1, 1, 4}, {i2, 
	1, 4}, {i3, 1, 4}, {i4, 1, 4}, {j, 1, 4}];
Do[VInt[i1, i2, i3, i4] = 
Sum[Frigus[(VV[i1, i2, i3, i4, j] /. 
x[j] -> 1) - (VV[i1, i2, i3, i4, j] /. x[j] -> -1)], {j, 1, 
	4}], {i1, 1, 4}, {i2, 1, 4}, {i3, 1, 4}, {i4, 1, 4}];
Do[Print[a[i1, i2, i3, i4] + 
Simplify[
VInt[i1, i2, i3, i4] - 
VInt[i3, i4, i1, 
i2], {(3 \[Sigma][1] + 6 \[Sigma][2] + 
	6 \[Sigma][3] + \[Sigma][4]) == 
	Pi^2/4, (3 \[Sigma][1] - 6 \[Sigma][2] - 
	6 \[Sigma][3] + \[Sigma][4]) == 0}]], {i1, 1, 4}, 
Do[CCC[i1, i2, i3, i4] = 
Simplify[
VInt[i1, i2, i3, i4] - 
VInt[i3, i4, i1, 
i2], {(3 \[Sigma][1] + 6 \[Sigma][2] + 
	6 \[Sigma][3] + \[Sigma][4]) == 
	Pi^2/4, (3 \[Sigma][1] - 6 \[Sigma][2] - 
	6 \[Sigma][3] + \[Sigma][4]) == 0}], {i1, 1, 4}, {i2, 1, 
	4}, {i3, 1, 4}, {i4, 1, 4}]; {i2, 1, 4}, {i3, 1, 4}, {i4, 1, 4}];		
	\end{lstlisting}
	\section{The list of constants, MMA ver.}
	In the following codes,
	\[\text{CVol}=-\frac{\mathcal{V}}{2\pi^2}.\]
	\begin{lstlisting}
Con[1,1,1,1]->CVol,Con[1,1,1,2]->0,Con[1,1,1,3]->0,
Con[1,1,1,4]->0,Con[1,1,2,1]->0,
Con[1,1,2,2]->
CVol-\[Zeta]11/9+\[Zeta]22/18+\[Zeta]33/18-\[Xi]11/
9+\[Xi]22/18+\[Xi]33/18,
Con[1,1,2,3]->-(\[Zeta]12/18)-\[Xi]12/18,
Con[1,1,2,4]->-(\[Zeta]13/18)-\[Xi]13/18,Con[1,1,3,1]->0,
Con[1,1,3,2]->-(\[Zeta]12/18)-\[Xi]12/18,
Con[1,1,3,3]->
CVol+\[Zeta]11/18-\[Zeta]22/9+\[Zeta]33/18+\[Xi]11/
18-\[Xi]22/9+\[Xi]33/18,
Con[1,1,3,4]->-(\[Zeta]23/18)-\[Xi]23/18,Con[1,1,4,1]->0,
Con[1,1,4,2]->-(\[Zeta]13/18)-\[Xi]13/18,
Con[1,1,4,3]->-(\[Zeta]23/18)-\[Xi]23/18,
Con[1,1,4,4]->
CVol+\[Zeta]11/18+\[Zeta]22/18-\[Zeta]33/9+\[Xi]11/
18+\[Xi]22/18-\[Xi]33/9,Con[1,2,1,1]->0,
Con[1,2,2,1]->\[Zeta]11/9-\[Zeta]22/18-\[Zeta]33/
18+\[Xi]11/9-\[Xi]22/18-\[Xi]33/18-Con[1,2,1,2],
Con[1,2,2,2]->0,
Con[1,2,3,1]->\[Zeta]12/18+\[Xi]12/18-Con[1,2,1,3],
Con[1,2,3,2]->-(\[Zeta]13/18)+\[Xi]13/18-Con[1,2,2,3],
Con[1,2,3,3]->-(\[Zeta]23/18)+\[Xi]23/18,
Con[1,2,4,1]->\[Zeta]13/18+\[Xi]13/18-Con[1,2,1,4],
Con[1,2,4,2]->\[Zeta]12/18-\[Xi]12/18-Con[1,2,2,4],
Con[1,2,4,3]->\[Zeta]22/6-\[Zeta]33/6-\[Xi]22/6+\[Xi]33/6-
Con[1,2,3,4],Con[1,2,4,4]->\[Zeta]23/18-\[Xi]23/18,
Con[1,3,1,1]->0,Con[1,3,1,2]->Con[1,2,1,3],
Con[1,3,2,1]->\[Zeta]12/18+\[Xi]12/18-Con[1,2,1,3],
Con[1,3,2,2]->\[Zeta]13/18-\[Xi]13/18,
Con[1,3,3,1]->-(\[Zeta]11/18)+\[Zeta]22/9-\[Zeta]33/
18-\[Xi]11/18+\[Xi]22/9-\[Xi]33/18-Con[1,3,1,3],
Con[1,3,3,2]->\[Zeta]23/18-\[Xi]23/18-Con[1,3,2,3],
Con[1,3,3,3]->0,
Con[1,3,4,1]->\[Zeta]23/18+\[Xi]23/18-Con[1,3,1,4],
Con[1,3,4,2]->-(\[Zeta]11/6)+\[Zeta]33/6+\[Xi]11/6-\[Xi]33/
6-Con[1,3,2,4],
Con[1,3,4,3]->-(\[Zeta]12/18)+\[Xi]12/18-Con[1,3,3,4],
Con[1,3,4,4]->-(\[Zeta]13/18)+\[Xi]13/18,Con[1,4,1,1]->0,
Con[1,4,1,2]->Con[1,2,1,4],
Con[1,4,1,3]->Con[1,3,1,4],
Con[1,4,2,1]->\[Zeta]13/18+\[Xi]13/18-Con[1,2,1,4],
Con[1,4,2,2]->-(\[Zeta]12/18)+\[Xi]12/18,
Con[1,4,3,1]->\[Zeta]23/18+\[Xi]23/18-Con[1,3,1,4],
Con[1,4,3,2]->\[Zeta]11/6-\[Zeta]22/6-\[Xi]11/6+\[Xi]22/6-
Con[1,4,2,3],Con[1,4,3,3]->\[Zeta]12/18-\[Xi]12/18,
Con[1,4,4,1]->-(\[Zeta]11/18)-\[Zeta]22/18+\[Zeta]33/
9-\[Xi]11/18-\[Xi]22/18+\[Xi]33/9-Con[1,4,1,4],
Con[1,4,4,2]->-(\[Zeta]23/18)+\[Xi]23/18-Con[1,4,2,4],
Con[1,4,4,3]->\[Zeta]13/18-\[Xi]13/18-Con[1,4,3,4],
Con[1,4,4,4]->0,Con[2,1,1,1]->0,
Con[2,1,1,2]->\[Zeta]11/9-\[Zeta]22/18-\[Zeta]33/
18+\[Xi]11/9-\[Xi]22/18-\[Xi]33/18-Con[1,2,1,2],
Con[2,1,1,3]->\[Zeta]12/18+\[Xi]12/18-Con[1,2,1,3],
Con[2,1,1,4]->\[Zeta]13/18+\[Xi]13/18-Con[1,2,1,4],
Con[2,1,2,1]->Con[1,2,1,2],Con[2,1,2,2]->0,
Con[2,1,2,3]->-(\[Zeta]13/18)+\[Xi]13/18-Con[1,2,2,3],
Con[2,1,2,4]->\[Zeta]12/18-\[Xi]12/18-Con[1,2,2,4],
Con[2,1,3,1]->Con[1,2,1,3],
Con[2,1,3,2]->Con[1,2,2,3],
Con[2,1,3,3]->-(\[Zeta]23/18)+\[Xi]23/18,
Con[2,1,3,4]->\[Zeta]22/6-\[Zeta]33/6-\[Xi]22/6+\[Xi]33/6-
Con[1,2,3,4],Con[2,1,4,1]->Con[1,2,1,4],
Con[2,1,4,2]->Con[1,2,2,4],
Con[2,1,4,3]->Con[1,2,3,4],
Con[2,1,4,4]->\[Zeta]23/18-\[Xi]23/18,
Con[2,2,1,1]->
CVol-\[Zeta]11/9+\[Zeta]22/18+\[Zeta]33/18-\[Xi]11/
9+\[Xi]22/18+\[Xi]33/18,Con[2,2,1,2]->0,
Con[2,2,1,3]->\[Zeta]13/18-\[Xi]13/18,
Con[2,2,1,4]->-(\[Zeta]12/18)+\[Xi]12/18,Con[2,2,2,1]->0,
Con[2,2,2,2]->CVol,Con[2,2,2,3]->0,Con[2,2,2,4]->0,
Con[2,2,3,1]->\[Zeta]13/18-\[Xi]13/18,Con[2,2,3,2]->0,
Con[2,2,3,3]->
CVol+\[Zeta]11/18+\[Zeta]22/18-\[Zeta]33/9+\[Xi]11/
18+\[Xi]22/18-\[Xi]33/9,
Con[2,2,3,4]->\[Zeta]23/18+\[Xi]23/18,
Con[2,2,4,1]->-(\[Zeta]12/18)+\[Xi]12/18,Con[2,2,4,2]->0,
Con[2,2,4,3]->\[Zeta]23/18+\[Xi]23/18,
Con[2,2,4,4]->
CVol+\[Zeta]11/18-\[Zeta]22/9+\[Zeta]33/18+\[Xi]11/
18-\[Xi]22/9+\[Xi]33/18,
Con[2,3,1,1]->-(\[Zeta]12/18)-\[Xi]12/18,
Con[2,3,1,2]->Con[1,2,2,3],
Con[2,3,1,3]->Con[1,3,2,3],
Con[2,3,1,4]->Con[1,4,2,3],
Con[2,3,2,1]->-(\[Zeta]13/18)+\[Xi]13/18-Con[1,2,2,3],
Con[2,3,2,2]->0,
Con[2,3,3,1]->\[Zeta]23/18-\[Xi]23/18-Con[1,3,2,3],
Con[2,3,3,2]->-(\[Zeta]11/18)-\[Zeta]22/18+\[Zeta]33/
9-\[Xi]11/18-\[Xi]22/18+\[Xi]33/9-Con[2,3,2,3],
Con[2,3,3,3]->0,
Con[2,3,4,1]->\[Zeta]11/6-\[Zeta]22/6-\[Xi]11/6+\[Xi]22/6-
Con[1,4,2,3],
Con[2,3,4,2]->-(\[Zeta]23/18)-\[Xi]23/18-Con[2,3,2,4],
Con[2,3,4,3]->-(\[Zeta]13/18)-\[Xi]13/18-Con[2,3,3,4],
Con[2,3,4,4]->\[Zeta]12/18+\[Xi]12/18,
Con[2,4,1,1]->-(\[Zeta]13/18)-\[Xi]13/18,
Con[2,4,1,2]->Con[1,2,2,4],
Con[2,4,1,3]->Con[1,3,2,4],
Con[2,4,1,4]->Con[1,4,2,4],
Con[2,4,2,1]->\[Zeta]12/18-\[Xi]12/18-Con[1,2,2,4],
Con[2,4,2,2]->0,Con[2,4,2,3]->Con[2,3,2,4],
Con[2,4,3,1]->-(\[Zeta]11/6)+\[Zeta]33/6+\[Xi]11/6-\[Xi]33/
6-Con[1,3,2,4],
Con[2,4,3,2]->-(\[Zeta]23/18)-\[Xi]23/18-Con[2,3,2,4],
Con[2,4,3,3]->\[Zeta]13/18+\[Xi]13/18,
Con[2,4,4,1]->-(\[Zeta]23/18)+\[Xi]23/18-Con[1,4,2,4],
Con[2,4,4,2]->-(\[Zeta]11/18)+\[Zeta]22/9-\[Zeta]33/
18-\[Xi]11/18+\[Xi]22/9-\[Xi]33/18-Con[2,4,2,4],
Con[2,4,4,3]->-(\[Zeta]12/18)-\[Xi]12/18-Con[2,4,3,4],
Con[2,4,4,4]->0,Con[3,1,1,1]->0,
Con[3,1,1,2]->\[Zeta]12/18+\[Xi]12/18-Con[1,2,1,3],
Con[3,1,1,3]->-(\[Zeta]11/18)+\[Zeta]22/9-\[Zeta]33/
18-\[Xi]11/18+\[Xi]22/9-\[Xi]33/18-Con[1,3,1,3],
Con[3,1,1,4]->\[Zeta]23/18+\[Xi]23/18-Con[1,3,1,4],
Con[3,1,2,1]->Con[1,2,1,3],
Con[3,1,2,2]->\[Zeta]13/18-\[Xi]13/18,
Con[3,1,2,3]->\[Zeta]23/18-\[Xi]23/18-Con[1,3,2,3],
Con[3,1,2,4]->-(\[Zeta]11/6)+\[Zeta]33/6+\[Xi]11/6-\[Xi]33/
6-Con[1,3,2,4],Con[3,1,3,1]->Con[1,3,1,3],
Con[3,1,3,2]->Con[1,3,2,3],Con[3,1,3,3]->0,
Con[3,1,3,4]->-(\[Zeta]12/18)+\[Xi]12/18-Con[1,3,3,4],
Con[3,1,4,1]->Con[1,3,1,4],
Con[3,1,4,2]->Con[1,3,2,4],
Con[3,1,4,3]->Con[1,3,3,4],
Con[3,1,4,4]->-(\[Zeta]13/18)+\[Xi]13/18,
Con[3,2,1,1]->-(\[Zeta]12/18)-\[Xi]12/18,
Con[3,2,1,2]->-(\[Zeta]13/18)+\[Xi]13/18-Con[1,2,2,3],
Con[3,2,1,3]->\[Zeta]23/18-\[Xi]23/18-Con[1,3,2,3],
Con[3,2,1,4]->\[Zeta]11/6-\[Zeta]22/6-\[Xi]11/6+\[Xi]22/6-
Con[1,4,2,3],Con[3,2,2,1]->Con[1,2,2,3],
Con[3,2,2,2]->0,
Con[3,2,2,3]->-(\[Zeta]11/18)-\[Zeta]22/18+\[Zeta]33/
9-\[Xi]11/18-\[Xi]22/18+\[Xi]33/9-Con[2,3,2,3],
Con[3,2,2,4]->-(\[Zeta]23/18)-\[Xi]23/18-Con[2,3,2,4],
Con[3,2,3,1]->Con[1,3,2,3],
Con[3,2,3,2]->Con[2,3,2,3],Con[3,2,3,3]->0,
Con[3,2,3,4]->-(\[Zeta]13/18)-\[Xi]13/18-Con[2,3,3,4],
Con[3,2,4,1]->Con[1,4,2,3],
Con[3,2,4,2]->Con[2,3,2,4],
Con[3,2,4,3]->Con[2,3,3,4],
Con[3,2,4,4]->\[Zeta]12/18+\[Xi]12/18,
Con[3,3,1,1]->
CVol+\[Zeta]11/18-\[Zeta]22/9+\[Zeta]33/18+\[Xi]11/
18-\[Xi]22/9+\[Xi]33/18,
Con[3,3,1,2]->-(\[Zeta]23/18)+\[Xi]23/18,Con[3,3,1,3]->0,
Con[3,3,1,4]->\[Zeta]12/18-\[Xi]12/18,
Con[3,3,2,1]->-(\[Zeta]23/18)+\[Xi]23/18,
Con[3,3,2,2]->
CVol+\[Zeta]11/18+\[Zeta]22/18-\[Zeta]33/9+\[Xi]11/
18+\[Xi]22/18-\[Xi]33/9,Con[3,3,2,3]->0,
Con[3,3,2,4]->\[Zeta]13/18+\[Xi]13/18,Con[3,3,3,1]->0,
Con[3,3,3,2]->0,Con[3,3,3,3]->CVol,Con[3,3,3,4]->0,
Con[3,3,4,1]->\[Zeta]12/18-\[Xi]12/18,
Con[3,3,4,2]->\[Zeta]13/18+\[Xi]13/18,Con[3,3,4,3]->0,
Con[3,3,4,4]->
CVol-\[Zeta]11/9+\[Zeta]22/18+\[Zeta]33/18-\[Xi]11/
9+\[Xi]22/18+\[Xi]33/18,
Con[3,4,1,1]->-(\[Zeta]23/18)-\[Xi]23/18,
Con[3,4,1,2]->Con[1,2,3,4],
Con[3,4,1,3]->Con[1,3,3,4],
Con[3,4,1,4]->Con[1,4,3,4],
Con[3,4,2,1]->\[Zeta]22/6-\[Zeta]33/6-\[Xi]22/6+\[Xi]33/6-
Con[1,2,3,4],Con[3,4,2,2]->\[Zeta]23/18+\[Xi]23/18,
Con[3,4,2,3]->Con[2,3,3,4],
Con[3,4,2,4]->Con[2,4,3,4],
Con[3,4,3,1]->-(\[Zeta]12/18)+\[Xi]12/18-Con[1,3,3,4],
Con[3,4,3,2]->-(\[Zeta]13/18)-\[Xi]13/18-Con[2,3,3,4],
Con[3,4,3,3]->0,
Con[3,4,4,1]->\[Zeta]13/18-\[Xi]13/18-Con[1,4,3,4],
Con[3,4,4,2]->-(\[Zeta]12/18)-\[Xi]12/18-Con[2,4,3,4],
Con[3,4,4,3]->\[Zeta]11/9-\[Zeta]22/18-\[Zeta]33/
18+\[Xi]11/9-\[Xi]22/18-\[Xi]33/18-Con[3,4,3,4],
Con[3,4,4,4]->0,Con[4,1,1,1]->0,
Con[4,1,1,2]->\[Zeta]13/18+\[Xi]13/18-Con[1,2,1,4],
Con[4,1,1,3]->\[Zeta]23/18+\[Xi]23/18-Con[1,3,1,4],
Con[4,1,1,4]->-(\[Zeta]11/18)-\[Zeta]22/18+\[Zeta]33/
9-\[Xi]11/18-\[Xi]22/18+\[Xi]33/9-Con[1,4,1,4],
Con[4,1,2,1]->Con[1,2,1,4],
Con[4,1,2,2]->-(\[Zeta]12/18)+\[Xi]12/18,
Con[4,1,2,3]->\[Zeta]11/6-\[Zeta]22/6-\[Xi]11/6+\[Xi]22/6-
Con[1,4,2,3],
Con[4,1,2,4]->-(\[Zeta]23/18)+\[Xi]23/18-Con[1,4,2,4],
Con[4,1,3,1]->Con[1,3,1,4],
Con[4,1,3,2]->Con[1,4,2,3],
Con[4,1,3,3]->\[Zeta]12/18-\[Xi]12/18,
Con[4,1,3,4]->\[Zeta]13/18-\[Xi]13/18-Con[1,4,3,4],
Con[4,1,4,1]->Con[1,4,1,4],
Con[4,1,4,2]->Con[1,4,2,4],
Con[4,1,4,3]->Con[1,4,3,4],Con[4,1,4,4]->0,
Con[4,2,1,1]->-(\[Zeta]13/18)-\[Xi]13/18,
Con[4,2,1,2]->\[Zeta]12/18-\[Xi]12/18-Con[1,2,2,4],
Con[4,2,1,3]->-(\[Zeta]11/6)+\[Zeta]33/6+\[Xi]11/6-\[Xi]33/
6-Con[1,3,2,4],
Con[4,2,1,4]->-(\[Zeta]23/18)+\[Xi]23/18-Con[1,4,2,4],
Con[4,2,2,1]->Con[1,2,2,4],Con[4,2,2,2]->0,
Con[4,2,2,3]->-(\[Zeta]23/18)-\[Xi]23/18-Con[2,3,2,4],
Con[4,2,2,4]->-(\[Zeta]11/18)+\[Zeta]22/9-\[Zeta]33/
18-\[Xi]11/18+\[Xi]22/9-\[Xi]33/18-Con[2,4,2,4],
Con[4,2,3,1]->Con[1,3,2,4],
Con[4,2,3,2]->Con[2,3,2,4],
Con[4,2,3,3]->\[Zeta]13/18+\[Xi]13/18,
Con[4,2,3,4]->-(\[Zeta]12/18)-\[Xi]12/18-Con[2,4,3,4],
Con[4,2,4,1]->Con[1,4,2,4],
Con[4,2,4,2]->Con[2,4,2,4],
Con[4,2,4,3]->Con[2,4,3,4],Con[4,2,4,4]->0,
Con[4,3,1,1]->-(\[Zeta]23/18)-\[Xi]23/18,
Con[4,3,1,2]->\[Zeta]22/6-\[Zeta]33/6-\[Xi]22/6+\[Xi]33/6-
Con[1,2,3,4],
Con[4,3,1,3]->-(\[Zeta]12/18)+\[Xi]12/18-Con[1,3,3,4],
Con[4,3,1,4]->\[Zeta]13/18-\[Xi]13/18-Con[1,4,3,4],
Con[4,3,2,1]->Con[1,2,3,4],
Con[4,3,2,2]->\[Zeta]23/18+\[Xi]23/18,
Con[4,3,2,3]->-(\[Zeta]13/18)-\[Xi]13/18-Con[2,3,3,4],
Con[4,3,2,4]->-(\[Zeta]12/18)-\[Xi]12/18-Con[2,4,3,4],
Con[4,3,3,1]->Con[1,3,3,4],
Con[4,3,3,2]->Con[2,3,3,4],Con[4,3,3,3]->0,
Con[4,3,3,4]->\[Zeta]11/9-\[Zeta]22/18-\[Zeta]33/
18+\[Xi]11/9-\[Xi]22/18-\[Xi]33/18-Con[3,4,3,4],
Con[4,3,4,1]->Con[1,4,3,4],
Con[4,3,4,2]->Con[2,4,3,4],
Con[4,3,4,3]->Con[3,4,3,4],Con[4,3,4,4]->0,
Con[4,4,1,1]->
CVol+\[Zeta]11/18+\[Zeta]22/18-\[Zeta]33/9+\[Xi]11/
18+\[Xi]22/18-\[Xi]33/9,
Con[4,4,1,2]->\[Zeta]23/18-\[Xi]23/18,
Con[4,4,1,3]->-(\[Zeta]13/18)+\[Xi]13/18,Con[4,4,1,4]->0,
Con[4,4,2,1]->\[Zeta]23/18-\[Xi]23/18,
Con[4,4,2,2]->
CVol+\[Zeta]11/18-\[Zeta]22/9+\[Zeta]33/18+\[Xi]11/
18-\[Xi]22/9+\[Xi]33/18,
Con[4,4,2,3]->\[Zeta]12/18+\[Xi]12/18,Con[4,4,2,4]->0,
Con[4,4,3,1]->-(\[Zeta]13/18)+\[Xi]13/18,
Con[4,4,3,2]->\[Zeta]12/18+\[Xi]12/18,
Con[4,4,3,3]->
CVol-\[Zeta]11/9+\[Zeta]22/18+\[Zeta]33/18-\[Xi]11/
9+\[Xi]22/18+\[Xi]33/18,Con[4,4,3,4]->0,
Con[4,4,4,1]->0,Con[4,4,4,2]->0,Con[4,4,4,3]->0,
Con[4,4,4,4]->CVol
	\end{lstlisting}
\section{The calculation about the $\sigma_i$}
	\begin{lstlisting}
J={{0,0,0,0},{0,0,0,0},{0,0,0,0},{0,0,0,0}};
x[1]=x1;
x[2]=x2;
x[3]=x3;
x[4]=x4;
\[Theta][1]=\[Theta]1;
\[Theta][2]=\[Theta]2;
\[Theta][3]=\[Theta]3;
x1=Cos[\[Theta]1];
x2=Sin[\[Theta]1]Cos[\[Theta]2];
x3=Sin[\[Theta]1]Sin[\[Theta]2]Cos[\[Theta]3];
x4=Sin[\[Theta]1]Sin[\[Theta]2]Sin[\[Theta]3];
R=(x1^4+x2^4+x3^4+x4^4);
Do[J[[i,j]]=D[x[i],\[Theta][j-1]],{i,1,4},{j,2,4}];
Do[J[[i,1]]=x[i],{i,1,4}];
MatrixForm[J]
Q=R*Det[J];
Simplify[Q]
Integrate[Q,{\[Theta]1,0,Pi},{\[Theta]2,0,Pi},{\[Theta]3,0,2Pi}]
	\end{lstlisting}
	\bibliographystyle{amsalpha}
	
	\bibliography{ALE}
\end{document}